\documentclass[11pt, a4paper]{amsart}

\usepackage[english]{babel}
\usepackage{amsmath,enumerate, amsfonts, amssymb,amsthm}
\usepackage[dvips]{graphics}
\usepackage{graphicx} 
\usepackage[dvipsnames]{xcolor}
\usepackage{hyperref}
\hypersetup{
	colorlinks=true,
	linkcolor={Blue},
	citecolor={red!60!black},
	urlcolor={cyan},
}
\usepackage{array}
\usepackage{float}
\usepackage{mathtools}
\usepackage{amsmath}
\usepackage{stmaryrd}

\makeatother
\makeatletter
\def\blfootnote{\xdef\@thefnmark{}\@footnotetext}
\makeatother

\newtheorem{thm}{Theorem}[section]
\newtheorem{lem}[thm]{Lemma}
\newtheorem{prop}[thm]{Proposition}
\newtheorem{cor}[thm]{Corollary}

\theoremstyle{definition}
\newtheorem{defi}[thm]{Definition}
\newtheorem{notation}[thm]{Notation}
\newtheorem{assumption}[thm]{Assumption}
\theoremstyle{remark} \newtheorem{ex}[thm]{Example}

\theoremstyle{remark} \newtheorem{rem}[thm]{Remark}

\title{Real log canonical thresholds and real log canonical thresholds with signs}

\author{Théo Jaudon}

\address{Université de Rennes, IRMAR - UMR 6625, F-35000 Rennes, France}

\email{theo.jaudon@univ-rennes.fr}

\begin{document}

\begin{abstract}
  In this paper, we study the properties of the real log canonical threshold $\text{rlct}(X,Y)$ associated with a pair $(X,Y)$ and the real log canonical thresholds with signs $\text{rlct}^{\pm}(f)$ associated with a polynomial function $f: \mathbb{R}^d \to \mathbb{R}$. In particular, we relate these invariants to other invariants arising from the theory of motivic integration. We also show how to compute these invariants in several particular cases, including the Newton non-degenerate case and the case of analytically irreducible plane curve germs.
\end{abstract}

\maketitle

\blfootnote{
\noindent\textbf{2020 Mathematics Subject Classification.} Primary: 14B05. Secondary: 14E18, 14P99.
}

\section*{Introduction}

Let $X$ be a smooth complex algebraic variety of dimension $d$ and let $Y, Z$ be algebraic subvarieties of $X$. The log canonical threshold of the pair $(X,Y)$ along $Z$ is defined by
$$\text{lct}_{Z}(X,Y) = \text{sup} \{ c \in \mathbb{Q}_+ \mid (X,c \cdot Y) \text{ is log canonical in a neighborhood of Z} \}.$$
It is an important invariant which appears in birational geometry in the study of the singularities of $Y$ in a neighborhood of $Z$ inside the ambient variety $X$, and it can be expressed in terms of a log resolution of the pair $(X,Y)$. One can also consider this invariant in the more general setting where $X$ is a complex analytic variety and $Y,Z$ are complex analytic subspaces of $X$. The case
$$\text{lct}_0(f) = \text{lct}_0(\mathbb{C}^d,\{f=0\})$$
where $f : (\mathbb{C}^d,0) \to (\mathbb{C},0)$ is a holomorphic function germ is of particular interest and admits an analytic characterization given by
$$\text{lct}_0(f) = \text{sup} \{ s> 0 \mid |f|^{-s} \text{ is locally } L^{2}  \text{ near }0\}.$$ This invariant arises naturally in a variety of contexts. For instance, it is shown in \cite{kollar1996singularitiespairs} that $-\text{lct}_0(f)$ is the largest root of the local Bernstein-Sato polynomial $b_{f,0}(s)$ associated with $f$. One can also mention the work of Mustaţă, Ein, de Fernex, Lazarsfeld, and others that relates the log canonical threshold to other invariants arising from the theory of motivic integration (see \cite{Mustata}, \cite{Mustata2}, \cite{Ein1}, \cite{Ein2}) as well as the connection between the ACC property for the set
$$\{ \text{lct}_0(f) \mid f:(\mathbb{C}^d,0) \to (\mathbb{C},0) \text{ holomorphic germ} \}$$
and the termination of flips in birational geometry (see \cite{DeFernex2}).
Asymptotic versions of the log canonical threshold, such as the $\alpha-$invariant, are also related to the existence of Kähler-Einstein metrics on Fano varieties (see \cite{Demailly-Kollar}) and can be used, for instance, to prove the non-rationality of certain Fano varieties as in \cite{DeFernex}. We also refer to \cite{mustata2011impangalecturenoteslog} for a survey on log canonical thresholds.

Now assume that $X$ is a smooth algebraic variety defined over $\mathbb{R}$ and that $Y$ and $Z$ are algebraic subvarieties of $X$ defined over $\mathbb{R}$. One can analogously consider the real log canonical threshold of the pair $(X,Y)$ along $Z$, which we denote by $\text{rlct}_Z(X,Y)$. This is an invariant of the singularities of the real locus of $Y$, that is, of $Y(\mathbb{R})$ endowed with its structure of real analytic space and this definition can in fact be generalized to the more general setting where $X$ is a real analytic variety and $Y,Z$ are real analytic subspaces of $X$. The case
$$\text{rlct}_0(f) = \text{rlct}_0(\mathbb{R}^d,\{f=0\})$$
where $f : (\mathbb{R}^d,0) \to (\mathbb{R},0)$ is a real analytic function germ is of particular interest and admits an analytic characterization given by
$$\text{rlct}_0(f) = \text{sup} \{ s> 0 \mid |f|^{-s} \text{ is locally } L^{1}  \text{ near }0 \text{ in } \mathbb{R}^d\}.$$
Although this invariant has received comparatively little attention in geometry (see however \cite{saito2007reallogcanonicalthresholds}), it plays, for instance, an important role in Bayesian statistics (see for example \cite{Watanabe1}, \cite{Watanabe2}, \cite{imai2019estimatingreallogcanonical}). The real setting also allows one to consider finer invariants, namely the real log canonical thresholds with signs associated with a real analytic function germ $f :(\mathbb{R}^d,0) \to (\mathbb{R},0)$, which are defined by
$$\text{rlct}^{+}_0(f) = \text{sup} \{ s > 0 \mid |f|^{-s} \text{ is locally } L^{1} \text{ near } 0 \text{ in }  \mathbb{R}^d \cap\{f \geq 0\}   \}$$
and similarly
$$\text{rlct}^{-}_0(f) = \text{sup} \{ s > 0 \mid |f|^{-s} \text{ is locally } L^{1} \text{ near } 0 \text{ in }  \mathbb{R}^d \cap\{f \leq 0\}    \}.$$

The aim of this paper is to study the properties of these real log canonical thresholds, highlighting the similarities and differences with the usual log canonical threshold, but also to see how these real invariants are related to other invariants arising in singularity theory, and more specifically in real algebraic and analytic geometry.

This paper is organized as follows. In the first section, we begin by introducing the definitions of the real log canonical threshold associated with a pair $(X,Y)$, where $X$ is a smooth algebraic variety over $\mathbb{R}$ and $Y$ is a closed subvariety of $X$. We also define the real log canonical thresholds with signs associated with an analytic function $f :\mathbb{R}^d \to \mathbb{R}$.\newline In the second section, we give several characterizations of these real log canonical thresholds in terms of invariants arising from the theory of motivic integration. For instance, we relate real log canonical thresholds with signs to the asymptotic behavior of the dimension of the real locus of certain contact spaces.

\begin{thm}[= Theorem \ref{thm 5}]
  One has a uniform bound in $n \in \mathbb{N} $ of the form
$$\text{dim } \mathcal{X}_n^{\pm}(\mathbb{R}) \leq nd+d-\text{rlct}^{\pm}(f)n -1$$
and $\text{rlct}^{\pm}(f)$ is the optimal exponent for uniform bounds as above. In other words, the real log canonical threshold is given by
    $$\text{rlct}^{\pm}(f)=\underset{n \geq 0}{\text{inf }}\frac{nd+d-\text{dim } \mathcal{X}_n^{\pm}(\mathbb{R}) -1}{n}.$$
  Moreover, the above inequality is in fact an equality for infinitely many integers, and hence one also has asymptotically
  $$\text{rlct}^{\pm}(f)=\underset{n \geq 0}{\text{ lim inf }}\frac{nd+d-\text{dim } \mathcal{X}_n^{\pm}(\mathbb{R}) -1}{n}.$$
\end{thm}

In the third section, we establish several properties satisfied by these real log canonical thresholds, and we also provide a number of examples illustrating the fact that real log canonical thresholds do not satisfy all of the more rigid properties of the usual log canonical threshold, such as inversion of adjunction and the lower semicontinuity property.\newline
In the fourth section, we describe a convenient way to read off the real log canonical threshold and the real log canonical thresholds with signs of a Nash function germ $f :(\mathbb{R}^d,0) \to (\mathbb{R},0)$ from its Newton polyhedron, under a non-degeneracy assumption. The results of this section allow us to provide a rich source of examples throughout the paper. In particular, we show how to extract the real log canonical thresholds from the Newton polyhedra. We refer the reader to Section \ref{section 4} for all the notation and definitions.
\begin{thm}[= Theorem \ref{thm 10}]
    Assume that $X_{\tau}^{c} = \emptyset$ for all $\tau \prec \Gamma_c(f)$. Then
$$\text{rlct}_0^{\pm}(f) = l_f^{\pm}.$$
Now assume that there exists at least one $\tau \prec \Gamma_c(f)$ such that $X_{\tau}^{c} \neq \emptyset$. Then
$$\text{rlct}_0^{\pm}(f) = \text{min}(1,l_f^{\pm}).$$
Moreover, the numbers $l_f^{\pm}$ can be easily read off the Newton polyhedron of $f$ by considering the intersection of the ray $\mathbb{R}_+(1,\dots,1)$ with the facets $\tau \prec \Gamma_c(f)$ for which $X_{\tau}^{\pm} \neq \emptyset$ (see \ref{prop 6} for details).
\end{thm}

In the final section, we show how to compute the pair $\{\text{rlct}_0^+(f),\text{rlct}_0^-(f)\}$ for an analytically irreducible function germ $f: (\mathbb{R}^2,0) \to (\mathbb{R},0)$. More precisely, let $(m;n_1,\ldots,n_g)$ denote the Puiseux characteristic exponents associated with $f$, where
 $m = \text{mult}_0(f)$. We define integers $\overline{n}_0,\dots,\overline{n}_g$ recursively by setting $\overline{n}_0 = m$, $\overline{n}_1 = n_1$, and for $j \geq 2$,
$$\overline{n}_j = n_j - n_{j-1} + \frac{e_{j-2}}{e_{j-1}}\overline{n}_{j-1}.$$ 
For $1 \leq i \leq g$, we also set $e_i =\gcd(m,n_1,\dots,n_i)$, and we adopt the conventions $e_0 = n_0 = m$.

\begin{thm}[= Theorem \ref{thm 8}]
     If $m$ is odd, then one has
$$\text{rlct}_0^{-}(f)=\text{rlct}_0^{+}(f)= \frac{1}{m} + \frac{1}{n_1}.$$

Now assume that $m$ is even and denote 
$$h = \text{min} \{ 1 \leq j \leq g \mid n_j \text{ is odd } \}.$$
Then 
$$\{ \text{rlct}_0^{+}(f),\text{rlct}_0^{-}(f) \} = \{ \frac{1}{m} + \frac{1}{n_1} , \frac{n_h+ m}{e_{h-1} \overline{n}_h} \}.$$

\end{thm}

\bigskip
\noindent
{\bf Acknowledgements.} This work is part of the author's PhD thesis. He thanks his advisor Goulwen Fichou for his suggestions and support during the preparation of this work. He would also like to thank Toshizumi Fukui for the insightful discussions on this topic.

\noindent \textbf{Data availability~} No data was used for the research described in this article.

\noindent \textbf{Conflict of interest~} The author has no conflict to disclose.

\noindent \textbf{Funding} The author have received support from the ANR New-Mirage ANR-23-CE40-0002-01.

\noindent \textbf{AI Disclosure statement} No AI tools were used during the preparation or writing of this manuscript.


\bigskip
\noindent

\section{Preliminaries}
\label{section 1}
In this section, we recall the basic definitions and notations related to real log canonical threshold (a standard reference is \cite{saito2007reallogcanonicalthresholds}) and we introduce the real log canonical thresholds with signs.
\\
\noindent
An algebraic variety over $\mathbb{R}$ will refer to a separated (not necessarily reduced) scheme of finite type over $\mathbb{R}$, while the term real algebraic variety will be reserved for varieties as defined in \cite{bcr}. If $X$ is an algebraic variety over $\mathbb{R}$, the set $X(\mathbb{R})$ of real closed points of $X$ is naturally endowed with a structure of real algebraic variety and, a fortiori, with a structure of real analytic space. Conversely, if $Y$ is a real algebraic variety, there exists an algebraic variety $X$ over $\mathbb{R}$ (generally not unique) such that $X(\mathbb{R})$ and $Y$ are isomorphic as real algebraic varieties. We denote by $\text{Var}_{\mathbb{R}}$ the category of algebraic varieties over $\mathbb{R}$ and by $\mathbb{R}\text{Var}$ the category of real algebraic varieties. Taking the real closed points gives a functor
\[
\begin{aligned}
\mathrm{Var}_{\mathbb{R}}
    &\longrightarrow \mathbb{R}\mathrm{Var}\\
X
    &\longmapsto X(\mathbb{R})
\end{aligned}
\]
but it does not induce an equivalence of categories, since $\mathbb{R}\text{Var}$ has strictly more morphisms than $\text{Var}_{\mathbb{R}}$. This larger class of morphisms implies, for example, that every quasi-projective real algebraic variety is affine (\cite{bcr} Theorem 3.4.4). Finally, if $X$ is an algebraic variety over $\mathbb{R}$, we denote by $X_{\mathbb{C}}=X \times_{\text{Spec}(\mathbb{R})} \text{Spec}(\mathbb{C})$ its complexification, and we simply write $X(\mathbb{C})$ for the set of complex points of $X_{\mathbb{C}}$.

\begin{defi}
    Let $X$ be a smooth real analytic manifold, $Y$ an analytic subspace of $X$ defined by an ideal sheaf $\mathcal{I}$, and let $x \in X$. The real log canonical threshold of the pair $(X,Y)$ at $x$ is defined by 
    $$\text{rlct}_x(X,Y) = \text{sup} \{ s > 0 \mid \frac{1}{\sum_{i=1}^p |f_i|^{s}} \text{ is locally integrable near } x   \} \in \mathbb{R}_+ \cup \{+ \infty \}$$
    where $f_1,\dots,f_p$ are generators of $\mathcal{I}_{x}$. One can check that this quantity does not depend on the choice of generators of $\mathcal{I}_x$. More generally, if $Z$ is an analytic subspace of $X$, the real log canonical threshold of the pair $(X,Y)$ along $Z$ is defined by
    $$\text{rlct}_Z(X,Y) = \underset{x \in Z}{\text{min}} \text{ rlct}_x(X,Y)$$
    and we simply write $\text{rlct}(X,Y) = \text{ rlct}_X(X,Y)$ for the global real log canonical threshold of the pair $(X,Y)$.
\end{defi}
This real log canonical threshold is a real analytic invariant that captures information about the singularities of $Y$ inside the ambient manifold $X$. A case that will often be of interest to us in what follows is that in which $X$ is a smooth algebraic variety of dimension $d$ defined over $\mathbb{R}$ and that $Y$ and $Z$ two closed subvarieties of $X$. We denote by $\text{rlct}_Z(X,Y)$ the real log canonical threshold of the pair $(X(\mathbb{R}),Y(\mathbb{R}))$ along $Z(\mathbb{R})$ where $X(\mathbb{R}), Y(\mathbb{R})$ and $Z(\mathbb{R})$ are naturally endowed with their structure of real analytic spaces. We also denote by $\text{lct}_Z(X,Y)$ the usual log canonical threshold of the pair $(X,Y)$ along $Z$, which in this case can also be characterized (see \cite{kollar1996singularitiespairs},\cite{Kollar_Mori_1998}) by  
$$\text{lct}_Z(X,Y) = \text{sup} \{ c \in \mathbb{Q}_+ \mid (X,c \cdot Y) \text{ is log canonical in a neighborhood of } Z   \} $$

A particular case that we will consider is the case when $Y$ is an algebraic hypersurface of $\mathbb{A}^d_{\mathbb{R}}$ defined by a polynomial function $f \in \mathbb{R}[x_1,\dots,x_d]$. In this situation, we simply write $\text{rlct}_Z(f)$ the real log canonical threshold of the pair $(\mathbb{A}^d_{\mathbb{R}},\{f=0\})$ along some variety $Z$. More generally, to a real analytic function germ $f :(\mathbb{R}^d,0) \to (\mathbb{R},0)$, one can associate the invariant $\text{rclt}_0(f)$.

\begin{notation}
    \label{not 1}
Let $X$ be a smooth algebraic variety defined over $\mathbb{R}$, and let $Y,Z$ be two closed subvarieties of $X$.
Let $ \sigma : (X',\sigma^{-1}(Y)) \to (X,Y)$ be a log resolution of $(X,Y)$, i.e., $X'$ is a smooth algebraic variety defined over $\mathbb{R}$ and $\sigma$ is a proper and birational map that is an isomorphism outside $Y$ such that the divisors $\sigma^{-1} (Y)$ and $K_{X'/X}=K_{X'} -\sigma^*K_X$ are simultaneously normal crossings.\\
Such a modification always exists according to Hironaka's theorem \cite{hiro1} and we will also say that $\sigma :(X^{'}, \sigma^{-1}(Y)) \to (X,Y)$ is an embedded resolution of $(X,Y)$.
Once we have a log resolution, we denote by $\underset{j \in J}{\bigcup} E_j$ the decomposition into irreducible components of $\sigma^{-1}(Y)$, where, by definition, the $E_j$ are smooth irreducible hypersurfaces in $X'$ that intersect transversally. We then denote by $\underset{I}{ \bigsqcup}E^{0}_I$ the canonical stratification of $\underset{j \in J}{\bigcup} E_j$ into smooth subvarieties, where $I$ runs over the set of non-empty subsets of $J$ and $E^{0}_I = \underset{i \in I}{\bigcap} E_i \setminus \underset{j \in J\setminus I}{\bigcup} E_j$. Finally, the numerical data of the resolution  are the integers $N_j$ and $\nu_j$ such that
$$\sigma^{-1}(Y) = \sum\limits_{j \in J} N_j E_j \text{ and } K_{X'/X} = \sum\limits_{j \in J} (\nu_j-1)E_j$$
and we will frequently consider the set
$$J_{\mathbb{R}}= \{ j \in J \mid E_j(\mathbb{R}) \neq \emptyset  \}.$$
\end{notation}

One can compute both the real and complex log canonical thresholds using a log resolution.

\begin{prop}
\label{prop 2}
    Let $\sigma : (X',\sigma^{-1}(Y)) \to (X,Y)$ be a log resolution of $(X,Y)$. Then 
    $$\text{lct}_Z(X,Y) = \text{min} \{ \frac{\nu_i}{N_i} \mid i \in J, E_i \cap \sigma^{-1}(Z) \neq \emptyset  \} \in \mathbb{Q}_+ \cup \{+\infty\}$$ 
    and 
    $$\text{rlct}_Z(X,Y) = \text{min} \{ \frac{\nu_i}{N_i} \mid i \in J, (E_i \cap \sigma^{-1}(Z))(\mathbb{R}) \neq \emptyset  \} \in \mathbb{Q}_+ \cup \{+\infty\}.$$ 
    
\end{prop}

\begin{proof}
These two equalities are well known and can be proved using the change of variables formula. We refer, for instance, to \cite{mustata2011impangalecturenoteslog} and \cite{saito2007reallogcanonicalthresholds}.
\end{proof}

\begin{ex}
\label{ex 1}
    If $Y$ is a smooth subvariety of $X$ whose real locus is non-empty, then 
    $$\text{rlct}(X,Y) = \text{lct}(X,Y) = d - \text{ dim } Y = d - \text{ dim } Y(\mathbb{R}).$$
    If $f(x_1,\dots,x_d) = x_1^{N_1}\dots x_d^{N_d}$ has global normal crossings, then $$\text{rlct}(f)=\text{lct}(f) = \underset{1 \leq i \leq d}{ \text{min}} \frac{1}{N_i}.$$
    If $f(x,y) = y^{2} - x^{3}$, one can compute a log resolution of $f$ to obtain 
    $\text{rlct}(f)=\text{lct}(f) = \frac{5}{6}.$\newline
    If $f(x,y) = x^{2} + y^{2} +1$ one finds that 
    $\text{rlct}(f)= + \infty \text{ while } \text{lct}(f) = 1.$\newline
    Take $f(x,y,z) = x^{2}+y^{2}+z^{2}$. One obtains a log resolution of $f$ by blowing up the origin and the exceptional divisor $E_1 \simeq \mathbb{P}^{2}_{\mathbb{R}}$ has numerical data $\nu_1 = 3$ and $N_1 = 2$ which yields $\text{lct}(f) = \min (1, \frac{3}{2}) = 1.$
    On the other hand, the real locus of $f$ is reduced to the origin, so that the strict transform of the zero locus of $f$ has no real points, and therefore $\text{rlct}(f) = \frac{3}{2}.$
\end{ex}

The following corollary follows immediately from Proposition \ref{prop 2}.

\begin{cor}
     One has the inequality
     $$\text{rlct}_Z(X,Y) \geq \text{lct}_Z(X,Y).$$
\end{cor}
However this inequality may be strict, as we have already seen in the previous example.\newline
To a real analytic function $f : \mathbb{R}^{d} \to \mathbb{R}$, one can associate both its usual log canonical threshold and its real log canonical threshold as discussed previously. Furthermore, one can also associate to $f$ real log canonical thresholds with signs in the following way.

\begin{defi}
    Let $f : \mathbb{R}^d \to \mathbb{R}$ be a real analytic function and $Z$ a real analytic subspace of $\mathbb{R}^d$. We define
   $$\text{rlct}^{+}_Z(f) = \text{sup} \{ s > 0 \mid |f|^{-s} \text{ is locally } L^{1} \text{ near } Z \text{ in }  \mathbb{R}^d \cap\{f \geq 0\}   \} \in \mathbb{R}_+ \cup \{ + \infty \}$$
    and similarly
$$\text{rlct}^{-}_Z(f) = \text{sup} \{ s > 0 \mid |f|^{-s} \text{ is locally } L^{1} \text{ near } Z \text{ in }  \mathbb{R}^d \cap\{f \leq 0\}    \} \in \mathbb{R}_+ \cup \{ + \infty \}.$$
\end{defi}

By definition one has 
$$\text{rlct}_Z(f) = \text{min} (\text{rlct}_Z^{+}(f),\text{rlct}_Z^{-}(f)).$$

\begin{notation}
\label{not 2}
    Whenever we are given a log resolution $\sigma : X' \to \mathbb{A}^d_{\mathbb{R}}$ of $f$, we will denote by $P(f)$ and $N(f)$ the sets defined by
$$P(f) = \{ x \in X'(\mathbb{R}) \mid (f \circ \sigma)(x) > 0 \} \text{ and } N(f) = \{ x \in X'(\mathbb{R}) \mid (f \circ \sigma)(x) < 0 \}.$$
We will also often use the notations
    $$J_{\mathbb{R}}^{+} = \{ j \in J \mid E_{j}(\mathbb{R}) \cap \overline{P(f)} \neq \emptyset \} \text{ and } J_{\mathbb{R}}^{-} = \{ j \in J \mid E_{j}(\mathbb{R}) \cap \overline{N(f)} \neq \emptyset \}$$
\end{notation}

These real log canonical thresholds with signs can also be expressed in terms of a log resolution associated to $f$.

\begin{prop}
    Let $\sigma : X' \to \mathbb{A}^d_{\mathbb{R}}$ be a log resolution of $f$. Then with the notation of \ref{not 2} one has
    $$\text{rlct}_Z^{+}(f) = \text{min} \{ \frac{\nu_i}{N_i} \mid i \in J, \overline{P(f)} \cap E_i(\mathbb{R}) \cap \sigma^{-1}(Z)(\mathbb{R}) \neq \emptyset \} \in \mathbb{Q}_+ \cup \{+\infty\}$$
    and similarly 
 $$\text{rlct}_Z^{-}(f) = \text{min} \{ \frac{\nu_i}{N_i} \mid i \in J, \overline{N(f)} \cap E_i(\mathbb{R}) \cap \sigma^{-1}(Z)(\mathbb{R}) \neq \emptyset \} \in \mathbb{Q}_+ \cup \{+\infty\}.$$
    In the global case one has 
    $$\text{rlct}^{\pm}(f) = \text{min} \{ \frac{\nu_i}{N_i} \mid i \in J_{\mathbb{R}}^{\pm} \}.$$
\end{prop}

\begin{proof}
   Both equalities can be proved in a similar way, so we will only show how to prove the first one. Note also that it suffices to establish this equality in the case where $Z$ is a point $x \in \mathbb{R}^d$, so we will therefore make this assumption from now on. We then consider the log resolution induced on real points $\sigma : X'(\mathbb{R}) \to \mathbb{R}^d$, where $\sigma$ is a proper map and an isomorphism outside a subset of measure zero in $\mathbb{R}^d$. By the change of variables formula, $|f|^{-s}$ is locally $L^1$ near $x$ in $\mathbb{R}^d \cap \{ f \geq 0 \}$ if and only if $|f \circ \sigma|^{-s}|\text{jac}(\sigma)|$ is locally $L^1$ near $\sigma^{-1}(x)$ in $X'(\mathbb{R}) \cap \{f \circ \sigma \geq 0\}$. Now the result follows from the fact that at every real point $p$ of $\sigma^{-1}(x)$, there exists a system of local coordinates $(y_1,\dots,y_d)$ centered at $p$ in $X'(\mathbb{R})$ and a system of local coordinates $(x_1,\dots,x_d)$ centered at $x$ in $\mathbb{R}^d$ such that one has 
$$\text{jac}(\sigma) = v(y)y_1^{\nu_1 -1}\dots y_d^{\nu_d -1} d y_1 \wedge \dots \wedge d y_d$$
where $v$ is a regular function satisfying $v(p) \neq 0$ and 
$$f ( \sigma (y)) = u_k(y) y_1^{N_1}\dots y_d^{N_d}$$
where $u_k$ is a regular function satisfying $u_k(p) \neq 0$.
\end{proof}

For simplicity, throughout the paper we assume that the ambient variety $X$ is smooth. However, one can more generally consider the usual log canonical threshold and the real log canonical threshold of a pair $(X,Y)$, where we only assume that $X$ is normal and $\mathbb{Q}$-Gorenstein.

\section{Jet schemes and motivic zeta functions}
\label{section 2}
In this section, we assume that $X$ is a smooth algebraic variety of dimension $d$ defined over $\mathbb{R}$, that $Y$ is a closed subvariety of $X$ defined by an ideal sheaf $\mathcal{I}$, and that $Z$ is a closed subvariety of $X$.
We recall some classical constructions from motivic integration in the case where the ambient variety $X$ is smooth, and we refer to \cite{denefloeser} for further details.

\begin{defi}
    We denote by $Y_n$ the algebraic variety defined over $\mathbb{R}$ parametrizing the $n$-jets on $Y$. In other words, $Y_n$ is characterized by the property that for every $\mathbb{R}$-algebra $A$, there is a natural bijection 
    $$\text{Hom}(\text{Spec } A, Y_n) \simeq \text{Hom}(\text{Spec } A[t]/(t^{n+1}) , Y).$$
    The construction of $Y_n$ is local, and one can describe $Y_n$ explicitly when $Y$ is affine. These $n$-jet spaces are thus defined intrinsically and encode information on the singularities of $Y$. In what follows, we will be particularly interested in the real locus $Y_n(\mathbb{R})$ of these jet spaces.
\end{defi}

\begin{rem}
    For the smooth algebraic variety $X$, the jet schemes have particularly nice properties. The $X_n$ are smooth algebraic varieties of dimension $(n+1)d$ for all $n \geq 0$, and the truncation morphisms  $\pi_{n',n} : X_{n'} \to X_n$ are smooth with fibers $\mathbb{A}_\mathbb{R}^{(n'-n)d}$ for every $n' > n \geq 0$.
\end{rem}

\begin{defi}
    Take $Y= \text{Spec}(\mathbb{R}[x,y,z]/(x_1^{2}+y_1^{2}+z_1^{2})$. One finds that 
    $$Y_1 \simeq T~ Y \simeq \text{Spec}(\frac{\mathbb{R}[x_1,y_1,z_1,x_2,y_2,z_2]}{(x^{2}+y^{2}+z^{2},x_1x_2 +y_1y_2+z_1z_2)})$$
    which is a scheme of dimension $4$, but whose real locus is 
    $Y_1(\mathbb{R}) \simeq\mathbb{R}^{3}$.
\end{defi}

\begin{defi}
    For $n \in \mathbb{N}$, we denote by $\mathcal{X}_n$ the set of $n$-jets of $X$ having contact of order $n$ along $Y$, i.e. 
    $$\mathcal{X}_n=\{\gamma \in X_n \mid\text{ord}_Y(\gamma) = n \}$$
    where $\text{ord}_Y(\gamma) = \text{min} \text{ ord}_t(g \circ \gamma)$ as $g$ ranges over local sections of $\mathcal{I}$ in a neighborhood of $\gamma(0)$.
    
\end{defi}

\begin{defi}
       We denote by K$_0(\text{Var}_{\mathbb{R}})$ (resp. K$_0$($\mathbb{R}$Var)) the free abelian group generated by the isomorphism classes $[X]$ of algebraic varieties over $\mathbb{R}$ (resp. of real algebraic varieties) modulo all relations of the form $[X]-[Y]- [X\setminus Y]$ whenever $Y$ is a closed subvariety of $X$. We equip K$_0(\text{Var}_{\mathbb{R}})$ and K$_0$($\mathbb{R}$Var) with a ring structure by setting $[X]\times [Y] = [X \times Y]$ and we refer to K$_0(\text{Var}_{\mathbb{R}})$ (resp. K$_0$($\mathbb{R}$Var)) as the Grothendieck ring of algebraic varieties over $\mathbb{R}$ (resp. the Grothendieck ring of real algebraic varieties). Finally, we denote by $\mathcal{M}_{\mathbb{R}}$ the localization of K$_0$($\mathbb{R}$Var) with respect to the class of the affine line $\mathbb{L}=[\mathbb{R}]$.
\end{defi}

\begin{rem}
\label{rem 3}
    There is a natural ring morphism 
  \[
\begin{aligned}
\text{K}_0(\mathrm{Var}_{\mathbb{R}})
    &\longrightarrow \text{K}_0(\mathbb{R}\mathrm{Var}) \\
[X]
    &\longmapsto [X(\mathbb{R})]
\end{aligned}
\]
which, after localization, also induces a morphism

 \[ 
\begin{aligned}
\text{K}_0(\mathrm{Var}_{\mathbb{R}})[\mathbb{A}^1_{\mathbb{R}}]^{-1}
    &\longrightarrow \text{K}_0(\mathbb{R}\mathrm{Var})[\mathbb{R}]^{-1} = \mathcal{M}_{\mathbb{R}} .
\end{aligned}
\]

One can analogously define the Grothendieck ring
 $\text{K}_0(\mathrm{Var}_{\mathbb{C}})$ and there is a base change morphism

\[ 
\begin{aligned}
\text{K}_0(\mathrm{Var}_{\mathbb{R}})
    &\longrightarrow \text{K}_0(\mathrm{Var}_{\mathbb{C}}) \\
[X]
    &\longmapsto [X_{\mathbb{C}}]
\end{aligned}
\]
\end{rem}

\begin{defi}
\label{def 1}
    The real motivic zeta function associated to $Y$ is defined by 
    $$Z_{mot}(Y;T) = \sum\limits_{n \geq 0} [\mathcal{X}_n(\mathbb{R})] \mathbb{L}^{-nd}T^n \in \mathcal{M}_{\mathbb{R}}[[T]].$$
\end{defi}

\begin{rem}
    In the case where $Y$ is a hypersurface of $\mathbb{A}^d_{\mathbb{R}}$ defined by a polynomial function $f \in \mathbb{R}[x_1,\dots,x_d]$, this real motivic zeta function is an algebro-geometric analogue of the real archimedean zeta function given by
    $$s\in \mathbb{C} \mapsto Z(f,\phi;s) = \int_{\mathbb{R}^d}|f(x)|^{s} \phi(x) dx$$ where
    $\phi \in \mathcal{C}^\infty_c(\mathbb{R}^d)$, which has been extensively studied by Barlet, Mardhy, Jeddi (for instance, in \cite{Barlet}, \cite{Barlet1993}).
    As in the archimedean case, we will see that the real log canonical threshold corresponds (up to sign) to the largest pole of the real motivic zeta function.
\end{rem}

\begin{thm}[\cite{denefloeser}]
\label{thm 2}
Let $\sigma : (X',\sigma^{-1}(Y)) \to (X,Y)$ be a log resolution of $(X,Y)$. Using the notation from \ref{not 1}, one has 
$$Z_{mot}(Y;T) = \sum\limits_{ I \subset J_{\mathbb{R}}} (\mathbb{L}-1)^ {|I|} [E^{0}_I(\mathbb{R}) ] \underset{i \in I}{\prod}\frac{\mathbb{L}^{-\nu_i}T^{N_i}}{1-\mathbb{L}^{-\nu_i}T^{N_i}}$$
\end{thm}

\begin{rem}
 
More precisely, the rationality formula as proved in \cite{loes1} is an equality in the ring $\text{K}_0(\text{Var}_{\mathbb{R}})[\mathbb{A}^1_{\mathbb{R}}]^{-1}[[T]]$. However, the desired equality described in the previous theorem follows by applying the ring homomorphism $\text{K}_0(\text{Var}_{\mathbb{R}})[\mathbb{A}^1_{\mathbb{R}}]^{-1}[[T]] \to \text{K}_0(\mathbb{R}\text{Var})[\mathbb{R}]^{-1}[[T]]$ given by taking the real closed points, as already mentioned in Remark \ref{rem 3}.
The motivic zeta function of $f$ is therefore a rational function; more precisely, $Z_{mot,0}(f;T)$ belongs to $\mathcal{M}_{\mathbb{R}} [T][\frac{1}{1-\mathbb{L}^{-\nu_i}T^{N_i}}]_{i \in I}.$ In particular, the sequence of $[\mathcal{X}_n(f)]$ is determined by a finite amount of data.\\
    On the other hand, this also proves that the above expression does not depend on the chosen resolution, since $Z_{mot,0}(f;T)$ was defined intrinsically in \ref{def 1}.
\end{rem}

\begin{rem}
    One can also define real local zeta functions relative to a subvariety $Z \subset X$ by 
    $$Z_{mot,Z}(Y;T) = \sum\limits_{n \geq 0} [\mathcal{X}_{n,Z}(\mathbb{R})] \mathbb{L}^{-nd}T^n$$
    where we consider the following contact spaces
    $$\mathcal{X}_{n,Z} = \{\gamma \in X_n \mid\text{ord}_Y(\gamma) = n \text { and } \gamma(0) \in Z\}.$$
    For instance, to a polynomial function germ $f :(\mathbb{R}^{d},0) \to (\mathbb{R},0)$ one can associate local motivic zeta functions in a neighborhood of the origin, which are well-defined. Once we have a log resolution $\sigma : (X',\sigma^{-1}(Y)) \to (X,Y)$, we also have the following rationality formula:
    $$Z_{mot,Z}(Y;T) = \sum\limits_{ I \subset J_{\mathbb{R}}} (\mathbb{L}-1)^ {|I|} [(E^{0}_I \cap \sigma^{-1}(Z))(\mathbb{R}) ] \underset{i \in I}{\prod}\frac{\mathbb{L}^{-\nu_i}T^{N_i}}{1-\mathbb{L}^{-\nu_i}T^{N_i}}.$$
\end{rem}

In practice, one usually works with more concrete specializations of the zeta functions, obtained via motivic measures on real algebraic varieties, that is, via ring morphisms K$_0(\mathbb{R}\text{Var}) \to A$. The finest additive and multiplicative invariant of real algebraic varieties known to date is the virtual Poincaré polynomial.

\begin{thm}[\cite{mccrory1} Corollary 2.2]
    \label{thm 4}
    There exists a unique ring morphism $ \beta : \text{K}_0(\mathbb{R}\text{Var}) \to \mathbb{Z}[u]$ such that $\beta(X) = \sum\limits_i \text{dim}(\text{H}_i(X;\mathbb{Z}/2\mathbb{Z})) u^{i}$ when $X$ is smooth and compact. One can recover the Euler characteristic with compact support of $X$ by evaluating $\beta(X)$ at $u=-1$ and the degree of $\beta(X)$ is equal to $\dim X(\mathbb{R})$. Moreover, the leading coefficient of $\beta(X)$ is equal to the number of connected components of the real locus of a smooth compactification of $X$, in particular it is always strictly positive whenever $\beta(X) \neq 0$.
\end{thm}

\begin{rem}
\label{rem 1}
By applying the virtual Poincaré polynomial, one obtains the specialization 
$$Z_{\beta}(Y;T) = \sum\limits_{ I \subset J_{\mathbb{R}}} (u-1)^ {|I|} \beta(E^{0}_I (\mathbb{R})) \underset{i \in I}{\prod}\frac{u^{-\nu_i}T^{N_i}}{1-u^{-\nu_i}T^{N_i}} \in \mathbb{Z}[u,u^{-1}](T)$$
which was first studied in \cite{goulw1}.
Evaluating $T$ at $u^{-s}$ and taking the limit as $u \to1$ yields the specialization
$$Z_{top}(Y;s)= \sum\limits_{ I \subset J_{\mathbb{R}}} \mu(E^{0}_I (\mathbb{R})) \underset{i \in I}{\prod}\frac{1}{\nu_i+sN_i} \in \mathbb{Q}(s)$$
studied in \cite{jaudon2026polesrealmotiviczeta}. We refer to section 4 of \cite{rodrigues3} for a precise definition of poles of $Z_{mot}(Y;T)$, and to \cite{jaudon2026polesrealmotiviczeta} for a definition of poles of $Z_{\beta}(Y;T)$. We will just use the fact that every pole of $Z_{top}(Y;T)$ is a pole of $Z_{\beta}(Y;T)$, and hence, a pole of $Z_{mot}(Y;T)$ (see \cite{jaudon2026polesrealmotiviczeta} Remark 1.20).
\end{rem}

\begin{defi}
To a real polynomial function $f \in \mathbb{R}[x_1,\dots,x_d]$, one can also associate zeta functions with signs of the form
$$Z_{mot}^{\pm}(f;T) = \sum\limits_{n \geq 0} [\mathcal{X}_n^{\pm}(\mathbb{R})] \mathbb{L}^{-nd}T^n$$ 
where $\mathcal{X}_n ^{\pm}$ denotes the set of $n$-jets $\gamma$ of $\mathbb{A}^d_{\mathbb{R}}$ such that $(f\circ\gamma)(t) = \pm t^n  ~(\text{mod } t^{n+1})$.
\end{defi}

One also has a rationnality formula
$$Z_{mot}^{\pm}(f;T) = \sum\limits_{ \emptyset \neq I \subset J} (\mathbb{L}-1)^{|I|-1} [\widetilde{E}_I^{0,\pm}(\mathbb{R})] \underset{i \in I}{\prod}\frac{\mathbb{L}^{-\nu_i}T^{N_i}}{1-\mathbb{L}^{-\nu_i}T^{N_i}}$$
which involves coverings $\widetilde{E}_I^{0,\pm}$ of the strata $E_I^{0}$, depending on the sign of $f$ in a neighborhood of these strata. In particular, one checks that $\widetilde{E}_I^{0,\pm}(\mathbb{R}) = \emptyset$ whenever there exists $i \in I$ such that $i \notin J_{\mathbb{R}}^{\pm}$, therefore one can write

 $$Z_{mot}^{\pm}(f;T) = \sum\limits_{ \emptyset \neq I \subset J_{\mathbb{R}}^{\pm}} (\mathbb{L}-1)^{|I|-1} [\widetilde{E}_I^{0,\pm}(\mathbb{R})] \underset{i \in I}{\prod}\frac{\mathbb{L}^{-\nu_i}T^{N_i}}{1-\mathbb{L}^{-\nu_i}T^{N_i}}.$$

Moreover, here again, one can consider local zeta functions with signs in a neighborhood of a subvariety $Z \subset \mathbb{A}^d_{\mathbb{R}}$. We refer for instance to \cite{goulw1} for further details on these zeta functions with signs.

The rationality formula (Theorem \ref{thm 2}) yields a set of candidate poles of $Z_{mot}(Y;T)$, namely 
$$\{ -\frac{\nu_i}{N_i} \mid  i \in J_{\mathbb{R}}\}.$$
In particular, we see that $-\text{rlct}(X,Y)$ appears as the largest candidate pole of $Z_{mot}(Y;T)$.
However, determining which of these candidates are actual poles of these zeta functions (in real or complex geometry) is, in general, a very difficult problem (see \cite{veys3} for a survey).

\begin{thm}
    \label{thm 3}
     \begin{enumerate}
         \item The largest pole of $Z_{mot}(Y;T)$ is always $-\text{rlct}(X,Y)$.
         \item Let $f \in \mathbb{R}[x_1,\dots,x_d]$ be a real polynomial function. Then the largest pole of $Z_{mot}^{\pm}(f;T)$ is always $-\text{rlct}^{\pm}(f)$.
     \end{enumerate}
\end{thm}

\begin{proof}
    We will only prove (2), as (1) can be proved in a similar way.
    By remark \ref{rem 1}, it suffices to show that $-\text{rlct}^{\pm}(f)$ is a pole of the zeta function $Z_{\beta}^{\pm}(f;T)$ defined at the level of the virtual Poincaré polynomial. Assume that $m \in \mathbb{N}^*$ is the order of $-\text{rlct}^{\pm}(f)$ as a candidate pole of $Z_{\beta}^{\pm}(f;T)$ , that is, $m$ is the largest integer such that there exists $I \subset J_{\mathbb{R}}^{\pm}$ with $|I| = m$,  $\widetilde{E}_I^{0,\pm}(\mathbb{R}) \neq \emptyset$ and such that $\text{rlct}^{\pm}(f)= \frac{\nu_i}{N_i}$ for all $i \in I$. We then consider the contribution of $\widetilde{E}_I^{0,\pm}(\mathbb{R})$ to the residue of $Z_{\beta}^{\pm}(f;T)$ at $-\text{rlct}^{\pm}(f)$ for any $I$ like this. More precisely, we first truncate the function $Z_{\beta}^{\pm}(f,T)$ to keep only the term
    $$\sum\limits_{ I \subset K \subset J_{\mathbb{R}}^{\pm}}(u-1)^{|K|-1} \beta(\widetilde{E}_K^{0,\pm}(\mathbb{R})) \prod\limits_{k \in K} \frac{u^{-\nu_k}T^{N_k}}{1-u^{-\nu_k}T^{N_k}}=\frac{1}{u-1}\sum\limits_{ I \subset K \subset J_{\mathbb{R}}^{\pm}} \beta(\widetilde{E}_K^{0,\pm}(\mathbb{R})) \prod\limits_{k \in K} \frac{u-1}{u^{\nu_k}T^{-N_k} -1}.$$
    The contribution $\mathcal{R}_I^{\pm}$ of $\widetilde{E}_I^{0,\pm}(\mathbb{R})$ to the residue is then given by
    $$\underset{T \to u^{\frac{\nu_i}{N_i}}}{\text{lim}}(T-u^{\frac{\nu_i}{N_i}})^m \frac{1}{u-1}\sum\limits_{ I \subset K \subset J_{\mathbb{R}}^{\pm}} \beta(\widetilde{E}_K^{0,\pm}(\mathbb{R})) \prod\limits_{k \in K} \frac{u-1}{u^{\nu_k}T^{-N_k} -1} $$
    which yields
    $$ \mathcal{R}_I^{\pm} = \frac{-(1-u)^{m-1} }{(N_i u^{\frac{\nu_i}{N_i}})^m} \sum\limits_{ I \subset K \subset J_{\mathbb{R}}^{\pm}} \beta(\widetilde{E}_K^{0,\pm}(\mathbb{R})) \prod\limits_{k \in K \setminus I} \frac{u-1}{u^{\alpha_k}-1}$$
    where $\alpha_k = \nu_k-\frac{\nu_i}{N_i}N_k \in \mathbb{Q}$. In what follows, we omit the nonzero factor $\frac{-(1-u)^{m-1}}{u^{\frac{m\nu_i}{N_i}}}$ which is common to all contributions for the candidate pole $-\frac{\nu_i}{N_i}$. It remains to show that $\mathcal{R}_I^{\pm}$ is nonzero, and that contributions coming from different strata cannot cancel each other. By definition of $\text{rlct}^{\pm}(f)$, all $\alpha_k$ are positive. Moreover, they are nonzero otherwise $- \text{rlct}^{\pm}(f)$ would be a candidate pole of order at least m+1. We may view $\mathcal{R}_I^{\pm}$ as a function of the real variable $u \in \mathbb{R}^{*}_+$ and it suffices to show that $\mathcal{R}_I^{\pm}(u)$ is asymptotically equivalent to a nonzero term as $u \to + \infty$. Let $I \subset K \subset J_{\mathbb{R}}$ be such that $\beta(\widetilde{E}_K^{0,\pm}(\mathbb{R})) \neq 0$, which is for instance the case for $K=I$ by assumption. 
    By theorem \ref{thm 4}, one has 
    $$\beta(\widetilde{E}_K^{0,\pm}(\mathbb{R})) \prod\limits_{k \in K \setminus I} \frac{u-1}{u^{\alpha_k}-1} \underset{u \to + \infty}{\sim} C_Ku^{\gamma_K} $$
    for some $C_K \in \mathbb{N}^{*}$ and $\gamma_K \in \mathbb{Q}$. It follows that
    $$\mathcal{R}_I^{\pm}(u) \underset{u \to +\infty}{\sim} C u^\gamma$$ for some $C \in \mathbb{Q}^{*}_{+}$ and $\gamma \in \mathbb{Q}$, which ensures that $\mathcal{R}_I^{\pm} \neq 0$. Since the leading coefficient appearing in this asymptotic expansion is strictly positive, it also follows that contributions coming from different strata to the residue of $Z_{\beta}^{\pm}(f;T)$ at $ - \text{rlct}^{\pm}(f)$ cannot cancel each other
\end{proof}

The local versions of the two previous results also hold and can be proved in a similar way. We also note that point (1) of Theorem \ref{thm 3} has been very recently improved in \cite{budur2026remarkablesubsetpolesmotivic}, where the authors exhibit remarkable subsets of poles of the real and complex zeta functions, which always contain the corresponding log canonical thresholds. Let us mention that it is not known (both in the real and complex case, except for curves) whether the corresponding log canonical threshold is always a pole at the level of the topological zeta function.

\begin{rem}
    It may happen that the real and complex log canonical thresholds coincide but that the orders as poles of the real and complex zeta functions respectively differ. For example, if we consider $f=x^{2}+y^{2}$ we have $\text{rlct}_0(f)=\text{lct}_0(f)=1$ but at the level of topological zeta functions we find that 
    $$Z_{top,0}(f;s)=\frac{1}{s+1} \text{ while } Z_{top,0}^{\mathbb{C}}(f;s)=\frac{1}{(s+1)^{2}}$$
    where $Z_{top,0}^{\mathbb{C}}(f;s)$ denotes the usual topological zeta function, i.e., the Denef–Loeser topological zeta function (see \cite{loes2}), which is defined by
    $$Z_{top,0}^{\mathbb{C}}(f;s) = \sum\limits_{\emptyset \neq I \subset J} \chi((E^{0}_I \cap \sigma^{-1}(0))(\mathbb{C})) \underset{i \in I}{\prod}\frac{1}{\nu_i+sN_i}.$$
    Such phenomena are mentioned in \cite{jaudon2026polesrealmotiviczeta}
    and already in \cite{Barlet} in the archimedean case.
    \end{rem}

\begin{rem}
\label{rem 6}
    We have seen that, to a polynomial function $f : \mathbb{R}^d \to \mathbb{R}$, one can associate motivic zeta functions $Z_{mot}(f;T)$ and $Z_{mot}^{\pm}(f;T)$ with coefficients in the ring $\text{K}_0(\mathbb{R}\mathrm{Var})[\mathbb{R}]^{-1}$. More generally, consider a Nash function $f :\mathbb{R}^d \to \mathbb{R}$, that is an analytic function which is algebraic over $\mathbb{R}[x_1,\dots,x_d]$, or equivalently a $\mathcal{C}^{\infty}$ function which is semi-algebraic (see \cite{bcr}). One can also associate motivic zeta functions $Z_{mot}(f;T)$ and $Z_{mot}^{\pm}(f;T)$ to $f$, with coefficients in the Grothendieck ring $\text{K}_0(AS)$ of arc-symmetric sets (see, for instance, \cite{goulw1} or \cite{campesato}). The same rationality formulas hold in exactly the same way, so one can express these zeta functions in terms of a Nash modification $\sigma : X' \to \mathbb{R}^d$ such that $f \circ \sigma$ and $\operatorname{jac} \sigma$ are simultaneously normal crossings. Finally, the virtual Poincaré polynomial has been extended to arc-symmetric sets in \cite{goulw1}. In particular, Theorem \ref{thm 3}, whose proof is carried out at the level of the virtual Poincaré polynomial, extends in a similar way to Nash functions $f : \mathbb{R}^d \to \mathbb{R}$. 
\end{rem}    

One of the main motivations for studying these zeta functions in real geometry is that they provide invariants for the classification of Nash function germs under the blow-Nash equivalence.
More precisely, it is shown in \cite{goulw1} (Theorem 4.9) that the functions $Z_{\beta,0}(f;T)$ and $Z_{\beta,0}^{\pm}(f;T)$ are invariants with respect to the blow-Nash equivalence for Nash functions germs $f :(\mathbb{R}^d,0) \to (\mathbb{R},0)$.
We have seen that the real log canonical thresholds (resp. with signs) are always poles of the real zeta functions (resp. with signs), hence we immediately obtain the following corollary.

\begin{cor}
\label{cor 1}
    Let $f:(\mathbb{R}^d,0) \to (\mathbb{R},0)$ be a Nash function germ. Then $\text{rlct}_0(f)$ and $\text{rlct}^{\pm}_0(f)$ are invariants for the blow-Nash equivalence.
\end{cor}

\begin{prop}
    Let $f:(\mathbb{R}^d,0) \to (\mathbb{R},0)$ be a real analytic function germ. Then $\text{rlct}_0(f)$ and $\text{rlct}_0^{\pm}(f)$ are invariants under right bi-Lipschitz equivalence.
\end{prop}

\begin{proof}
    Let $h:(\mathbb{R}^d,0) \to (\mathbb{R}^d,0)$ be a bi-Lipschitz homeomorphism germ. Then, by Rademacher's theorem, $h$ is differentiable outside a set of measure zero, and there exists $K \in \mathbb{R}_+$ such that 
    $$|\text{jac}(h(x))| \leq K$$
    for all $x$ in a neighborhood of the origin, outside a subset of measure zero. For $s > 0$, the change of variables formula ensures that $|f|^{-s}$ is $L^1$ near the origin if and only if $|f \circ h |^{-s} |\text{jac} (h)|$ is $L^1$ near the origin, which is the case whenever $|f \circ h|^{-s}$ is $L^1$ near the origin, by the above inequality. Therefore $\text{rclt}_0(f \circ h) \leq \text{rlct}_0(f)$ and one has the equality $\text{rlct}_0(f) = \text{rlct}_0(f \circ h)$ by symmetry. Similarly, one shows that $\text{rlct}_0^{\pm}(f) = \text{rlct}_0^{\pm}(f \circ h)$.
\end{proof}

Using the rationality formula, one can also relate the real log canonical threshold with signs of $f$ to the dimension of the real locus of the contact spaces $\mathcal{X}_n^{\pm}$.

\begin{thm}
\label{thm 5}
One has a uniform bound in $n \in \mathbb{N} $ of the form
$$\text{dim } \mathcal{X}_n^{\pm}(\mathbb{R}) \leq nd+d-\text{rlct}^{\pm}(f)n -1$$
and $\text{rlct}^{\pm}(f)$ is the optimal exponent for uniform bounds as above. In other words, the real log canonical threshold is given by
    $$\text{rlct}^{\pm}(f)=\underset{n \geq 0}{\text{inf }}\frac{nd+d-\text{dim } \mathcal{X}_n^{\pm}(\mathbb{R}) -1}{n}.$$
  Moreover, the above inequality is in fact an equality for infinitely many integers, and hence one also has asymptotically
  $$\text{rlct}^{\pm}(f)=\underset{n \geq 0}{\text{ lim inf }}\frac{nd+d-\text{dim } \mathcal{X}_n^{\pm}(\mathbb{R}) -1}{n}.$$
\end{thm}

In order to prove Theorem \ref{thm 5}, we will temporarily denote by 
$$\alpha_0^{\pm}= \underset{n \geq 0}{\text{inf }}\frac{nd+d-\text{dim } \mathcal{X}_n^{\pm}(\mathbb{R})-1 }{n}$$
and we will prove below that $\alpha_0^{\pm}=\text{rlct}^{\pm}(f)$. For $\alpha \in \mathbb{Q}$ one can consider the generating series 
$$Z_{\beta}^{\pm}(f;u^{\alpha}T)u^{1-d} = \sum\limits_{n \geq 0} \beta(\mathcal{X}_n^{\pm}(\mathbb{R}))u^{-(nd-n\alpha +d-1)} T^n$$ 
whose coefficients are elements of $\bigcup_{k \geq 1}\mathbb{Z}[u^{\frac{1}{k}}, u^{-\frac{1}{k}}]$. In particular, these coefficients have a well-defined degree in $\mathbb{Q}$ and admit a limit in $\mathbb{R}\cup \{\pm \infty\}$ as $u \to +\infty$.
One can also characterize $\alpha_0^{\pm}$ as follows.

\begin{lem}
\label{lemme 1}
    One has 
    $$\alpha_0^{\pm} = \text{sup } \{\alpha \in \mathbb{Q} \mid \text{ Every coefficient of } Z_{\beta}^{\pm}(f;u^\alpha T)u^{1-d} \text{ tend to } 0 \text{ as } u \to +\infty\}, $$
    in other words
    $$\alpha_0^{\pm} =  \text{sup } \{\alpha \in \mathbb{Q} \mid \text{ Every coefficient of } Z_{\beta}^{\pm}(f;u^\alpha T)u^{1-d} \text{ has degree } < 0\}.$$
\end{lem}

\begin{proof}
    This follows from the definition of $\alpha_0^{\pm}$, using the fact that 
    $$Z_{\beta}^{\pm}(f;u^{\alpha}T)u^{-d} = \sum\limits_{n \geq 0} \beta(\mathcal{X}_n^{\pm}(\mathbb{R}))u^{-(nd-n\alpha +d-1)} T^n$$
   whose $n$-th coefficient has degree $\dim \mathcal{X}_n^{\pm}(\mathbb{R}) - (nd-n\alpha+d-1)$, since the degree of $\beta(\mathcal{X}_n^{\pm}(\mathbb{R}))$ is equal to $\dim \mathcal{X}_n^{\pm}(\mathbb{R})$ (see Theorem \ref{thm 4}).
\end{proof}

\begin{proof}[Proof of theorem \ref{thm 5}]
Using the rationality formula and the characterization in the previous lemma, we will show that $\alpha_0^{\pm} = \text{rlct}^{\pm}(f) = \underset{i \in J_{\mathbb{R}}^{\pm}}{\text{min }}\frac{\nu_i}{N_i}$. For $\alpha \in \mathbb{Q}$, the rationnality formula yields 
$$Z_{\beta}^{\pm}(Y;u^{\alpha}T)u^{1-d} = \sum\limits_{ I \subset J_{\mathbb{R}}^{\pm}} (u-1)^ {|I|}u^{-d} \beta(\widetilde{E}_I^{0,\pm}(\mathbb{R})) \underset{i \in I}{\prod}\frac{u^{-\nu_i + \alpha N_i}T^{N_i}}{1-u^{-\nu_i + \alpha N_i}T^{N_i}}.$$
Assume that $\alpha < \text{rlct}^{\pm}(f)$ and fix $I \subset J_{\mathbb{R}}^{\pm}$. The real locus of the smooth open stratum $E_I^0(\mathbb{R})$ has dimension $d-|I|$ and the same holds for the covering space $\widetilde{E}_I^{0,\pm}(\mathbb{R})$, hence $(u-1)^{|I|}u^{-d} \beta(\widetilde{E}_I^{0,\pm}(\mathbb{R}))$ has degree $0$. On the other hand, each coefficient of the generating series
$$\underset{i \in I}{\prod}\frac{u^{-\nu_i + \alpha N_i}T^{N_i}}{1-u^{-\nu_i + \alpha N_i}T^{N_i}}$$
is a finite sum of monomials of degree $\sum\limits_{k=1}^{|I|} m_k(-\nu_k + \alpha N_k) < 0$, where $m_1,\dots,m_{|I|} \in \mathbb{N}^*$. It follows that that all coefficients of
$$(u-1)^ {|I|}u^{-d} \beta(\widetilde{E}_I^{0,\pm}(\mathbb{R})) \underset{i \in I}{\prod}\frac{u^{-\nu_i + \alpha N_i}T^{N_i}}{1-u^{-\nu_i + \alpha N_i}T^{N_i}}$$
are of degree $< 0$. This is therefore also the case for $Z_{\beta}^{\pm}(f;u^{\alpha}T)u^{1-d}$, which is the sum of these terms over $I \subset J_{\mathbb{R}}^{\pm}$.\\
We now assume that $\alpha > \text{rlct}^{\pm}(f)$, and we need to exhibit a coefficient of $Z_{\beta}^{\pm}(f;u^{\alpha}T)u^{1-d}$ which does not tend to $0$ as $u \to +\infty$. As seen above, the rationality formula allows us to write 
$$Z_{\beta}^{\pm}(Y;u^{\alpha}T)u^{1-d} =\sum\limits_{n \geq 0}    \sum\limits_{I \subset J_{\mathbb{R}}^{\pm}} P_{n,I}^{\pm}(u) T^n $$ 
where $P_{n,I}^{\pm}(u)$ are elements of $\underset{k \geq 1}{\bigcup}\mathbb{Z}[u^{\frac{1}{k}}, u^{-\frac{1}{k}}]$ whose leading coefficient is strictly positive. It therefore suffices to exhibit $I \subset J_{\mathbb{R}}^{\pm}$ and $n \in \mathbb{N}$ such that $\text{deg}(P_{n,I}^{\pm}) > 0$. By assumption, there exists at least one component $E_1$ such that $\widetilde{E}_1^{0,\pm}(\mathbb{R})$ and $\alpha > \frac{\nu_1}{N_1}$. The term corresponding to $\widetilde{E}_1^{0,\pm}(\mathbb{R})$ is given by 
$$(u-1)u^{-d}\beta(\widetilde{E}_1^{0,\pm}(\mathbb{R}))\frac{u^{-\nu_1 + \alpha N_1}T^{N_1}}{1-u^{-\nu_1 + \alpha N_1}T^{N_1}} = \sum\limits_{k \geq 1} (u-1)u^{-d}\beta(\widetilde{E}_1^{0,\pm}(\mathbb{R})) u^{k(-\nu_1+\alpha N_1)} T^{kN_1}$$ and for instance $P_{N_1,1}^{\pm}(u)$, which is given by
$$P_{N_1,1}^{\pm}(u) = (u-1)u^{-d}\beta(\widetilde{E}_1^{0,\pm}(\mathbb{R})) u^{-\nu_1+\alpha N_1}$$
has strictly positive degree. With the same notation, we also see that the equality
$$\text{dim } \mathcal{X}_n^{\pm}(\mathbb{R}) = nd+d-\text{rlct}^{\pm}(f)n - 1$$
for example holds for all integers $n$ which are multiples of $N_1$.
\end{proof}

\begin{thm}
\label{thm 11}
One has a uniform bound in $n \in \mathbb{N} $ of the form
$$\text{dim } \mathcal{X}_n(\mathbb{R}) \leq nd+d-\text{rlct}(X,Y)n$$
and $\text{rlct}(X,Y)$ is the optimal exponent for uniform bounds as above. In other words, the real log canonical threshold is given by
    $$\text{rlct}(X,Y)=\underset{n \geq 0}{\text{inf }}\frac{nd+d-\text{dim } \mathcal{X}_n(\mathbb{R}) }{n}.$$
  Moreover, the above inequality is in fact an equality for infinitely many integers, and hence one also has asymptotically
  $$\text{rlct}(X,Y)=\underset{n \geq 0}{\text{ lim inf }}\frac{nd+d-\text{dim } \mathcal{X}_n(\mathbb{R}) }{n}.$$
\end{thm}

\begin{proof}
    The proof is the same as that of theorem \ref{thm 5} using the rationnality formula for naive zeta functions (Theorem \ref{thm 2}).
\end{proof}

\begin{rem}
    As already mentioned in Remark \ref{rem 6}, Theorems \ref{thm 5} and \ref{thm 11} also hold for Nash functions $f : \mathbb{R}^d \to \mathbb{R}$, and the proof is exactly the same, using the virtual Poincaré polynomial extended to arc-symmetric sets.
\end{rem}

Such results on the asymptotic behavior of the real loci of contact spaces already appear in \cite{goulw4} in the Newton non-degenerate case, where they are proved purely combinatorially. The authors use these results to show that the weights are blow-Nash invariants for quasi-homogeneous polynomials in three variables.

\begin{rem}
   In the local setting, one shows in a similar way that 
   $$\text{rlct}_Z^{\pm}(f)  =  \underset{n \geq 0}{\text{inf }}\frac{nd+d-\text{dim } \mathcal{X}_{n,Z}^{\pm}(\mathbb{R}) -1 }{n}  = \underset{n \geq 0}{\text{ lim inf }}\frac{nd+d-\text{dim } \mathcal{X}_{n,Z}^{\pm}(\mathbb{R})  -1}{n}$$
   and that
    $$\text{rlct}_Z(X,Y)  =  \underset{n \geq 0}{\text{inf }}\frac{nd+d-\text{dim } \mathcal{X}_{n,Z}(\mathbb{R}) }{n}  = \underset{n \geq 0}{\text{ lim inf }}\frac{nd+d-\text{dim } \mathcal{X}_{n,Z}(\mathbb{R}) }{n}.$$
    It may however happen that the inequality
    $$\text{dim } \mathcal{X}_{n,Z}(\mathbb{R}) \leq nd+d-\text{rlct}_Z(X,Y)n$$
    is never an equality and similarly in the case with signs. Let us consider, for example, $f(x,y)=x^3y^4$, so that
$\text{rlct}_0(f)=\frac{1}{4}$. Then, for $n \in \mathbb{N}^*$
such that $n \notin 3\mathbb{N}^*+4\mathbb{N}^*$, one has
$\mathcal{X}_{n,0}(\mathbb{R})=\emptyset$. On the other hand,
for $n \in \mathbb{N}^*$ that can be written as $n=3m+4l$, one finds that
    $$\dim(\mathcal{X}_{n,0}(\mathbb{R})) = 2+2n-m-l < 2 +2n - \frac{n}{4}.$$
\end{rem}

\begin{defi}
     We define the real Poincaré series associated with $Y$ by 
$$J(Y;T) = \sum\limits_{n \geq 0}[Y_n(\mathbb{R})]T^n \in \text{K}_0(\mathbb{R}\text{Var})[[T]].$$
\end{defi}
This series encodes information about the real loci of the $n$-jets spaces of $Y$ and is independent of the ambient variety $X$.
We then have the usual functional equation (cf. \cite{Mourtada2023})
$$J(Y;T)=\frac{Z_{mot}(Y;\mathbb{L}^dT)-[X]}{\mathbb{L}^dT-1}$$
which proves in particular the rationality of $J$. Using this functional equation together with Theorem \ref{thm 5}, one easily proves the following result, which is the real analogue of Corollary 3.4 in \cite{Mustata}.

\begin{prop}
    \label{prop 1}
    One has 
    $$\text{rlct}(X,Y)=d-\underset{n \geq 0}{\text{sup }}\frac{\text{dim } Y_n(\mathbb{R})}{n+1} = d-\underset{n \geq 0}{\text{lim sup }}\frac{\text{dim } Y_n(\mathbb{R})}{n+1} $$
Moreover, whenever $n$ is divisible by a finite set of integers depending on the resolution, one has the equality $\text{dim }Y_n(\mathbb{R}) =(n+1)(d-\text{rlct}(X,Y))$.   
\end{prop}

As is the case for the usual log canonical threshold, Proposition \ref{prop 1} in fact shows that $\text{rlct}(X,Y)$ depends only on $Y$ and on the ambient dimension $d = \text{dim } X$.

\begin{rem}
    Corollary 3.4 in \cite{Mustata} provides the uniform bound
    $$\text{dim } Y_n \leq (n+1)(d-\text{lct}(X,Y))$$
    where $\text{dim } Y_n$ is the scheme-theoretic dimension of $Y_n$, that is, the dimension of $Y_n$ viewed as a complex algebraic variety. Moreover, this bound is optimal in the sense that $\text{lct}(X,Y)$ is the largest possible coefficient for uniform bounds as above. On the other hand, if one restricts to the real locus of the jet spaces, the above bound still holds but is no longer necessarily optimal, as soon as $\text{rlct}(X,Y) > \text{lct}(X,Y)$. If we consider for example $Y \subset \mathbb{A}^{3}_{\mathbb{R}}$ defined by $x^{2}+y^{2}+z^{2}=0$ where $\text{lct}(\mathbb{A}^{3}_{\mathbb{R}},Y) =1$, one has an optimal uniform bound given by
    $$\text{dim }Y_n \leq 2(n+1).$$
    On the other hand, $\text{rlct}(\mathbb{A}^{3}_{\mathbb{R}},Y) =\frac{3}{2}$, which yields the optimal uniform bound
    $$\text{dim } Y_n(\mathbb{R}) \leq \frac{3}{2}(n+1).$$
   One can verify these estimates by computing $Z_{mot}(Y;T)$ and $J(Y;T)$ explicitly from a log-resolution of $Y$.
\end{rem}

\section{Properties and examples}
\label{section 3}
In this section, we state some properties of these real log canonical thresholds analogous to those of the complex log canonical threshold, which can be obtained directly from Proposition \ref{prop 1}, similarly to what was done in \cite{Mustata}. We also highlight the fact that these real log canonical thresholds do not, in general, satisfy some of the stronger properties of the usual log canonical threshold, such as lower semicontinuity and inversion of adjunction.

\begin{prop}
    One always has $\text{rlct}(X,Y) \leq d - \text{dim }Y(\mathbb{R})$.
\end{prop}

\begin{proof}
    By proposition \ref{prop 1}, one has 
    $$\text{rlct}(X,Y) = d-\underset{n \geq 0}{\text{sup }}\frac{\text{dim } Y_n(\mathbb{R})}{n+1}  \leq d - \text{dim }Y_0(\mathbb{R})$$
    where $Y_0(\mathbb{R}) \simeq Y(\mathbb{R})$.
\end{proof}

\begin{rem}
    For the complex log canonical threshold one has $\text{lct}(X,Y) \leq d - \text{dim } Y$, which in particular gives $\text{lct}(X,Y) \leq 1$ when $Y$ is a hypersurface of $X$, but this is not true in general in the real case. Such phenomena of real log canonical threshold greater than $1$ in the case where the real locus has codimension greater than $1$ is already mentioned in \cite{saito2007reallogcanonicalthresholds}.
    It may also happen that
    $$\text{lct}(X,Y) < \text{rlct}(X,Y) \leq 1.$$
    For example if we consider $f(x,y,z) = (x^{2}+y^{2}+z^{2})^k$ with $k \geq 2$, one finds that 
    $$\text{rlct}(f) =\frac{3}{2k} \text{ and } \text{lct}(f) =\frac{1}{k}.$$ 
\end{rem}

\begin{prop}
    Let $Y' \subset Y \subset X$ be closed subvarieties. Then one has
    $$\text{rlct}(X,Y') \geq \text{rlct}(X,Y).$$
\end{prop}

\begin{proof}
    The fact that $Y'$ is a closed subvariety of $Y$ implies that $Y_n'$ is a closed subvariety of $Y_n$ for all $n \in \mathbb{N}$, and it follows that 
    $$\text{rlct}(X,Y')= d-\underset{n \geq 0}{\text{sup }}\frac{\text{dim } Y'_n(\mathbb{R})}{n+1} \geq d-\underset{n \geq 0}{\text{sup }}\frac{\text{dim } Y_n(\mathbb{R})}{n+1} =\text{rlct}(X,Y). $$
\end{proof}

\begin{prop}
    Let $X,X'$ be smooth algebraic varieties defined over $\mathbb{R}$ of dimension $d$ and $d'$ respectively, and let $Y,Y'$ be closed subvarieties of $X$ and $X'$ respectively. Then
    $$\text{rlct}(X\times X',Y \times Y') = \text{rlct}(X,Y) + \text{rlct}(X',Y').$$
\end{prop}

\begin{proof}
    By proposition \ref{prop 1}, there exists an integer $n$ that simultaneously realizes the three equalities $\text{rlct}(X\times X',Y \times Y')= d+d' - \frac{\text{dim }( Y \times Y')_n(\mathbb{R})}{n+1}$, $\text{rlct}(X,Y) = d-\frac{\text{dim } Y_n(\mathbb{R})}{n+1}$ and $\text{rlct}(X',Y')=d'-\frac{\text{dim } Y_n'(\mathbb{R})}{n+1}$. Now, by functoriality, one has natural isomorphisms 
    $$(Y \times Y')_n(\mathbb{R}) \simeq Y_n(\mathbb{R}) \times Y'_n(\mathbb{R})$$ so the result follows.
\end{proof}

On the other hand, the real log canonical thresholds does not satisfy properties as strong as those of the usual log canonical threshold, as illustrated by the few examples and remarks below.

\begin{rem}
    The complex log canonical threshold satisfies a lower semicontinuity property (see \cite{Demailly-Kollar}, \cite{Mustata}) but this is not the case for the real log canonical threshold, and a fortiori not for the real log canonical thresholds with signs. For instance, consider the one-parameter family defined by $f_t(x,y,z) = x^{2}+y^{2}+tz^{2}$. Then one finds that $\text{rlct}(f_t) =1$ for $t \leq 0$ while $\text{rlct}(f_t) = \frac{3}{2}$ for $t > 0$. On the other hand, $\text{lct}(f_t) = 1$ for all $t$.
\end{rem}

\begin{prop}
    Let $I$ be a Nash manifold and $ F : (x,t) \in \mathbb{R}^d \times I \mapsto f_t(x) \in  \mathbb{R}$ be a polynomial function. Assume that for all $t \in I$, $f_t$ has an isolated critical point at $0$. Then there exists a stratification
    $I = \underset{l \in L}{\bigsqcup}S_l$ into Nash submanifolds such that the maps $t \mapsto\text{rlct}_0(f_t)$ et $t \mapsto\text{rlct}^{\pm}_0(f_t)$ are constant on each connected component of $S_l$.
\end{prop}

\begin{proof}
    Theorem VI of \cite{Fukui2} provides a stratification as above such that the family $(f_t)_{t \in I}$ is trivial for blow-Nash equivalence on each connected component of $S_l$, and the result then follows from Corollary \ref{cor 1}.
\end{proof}

\begin{rem}
    If $H \subset X$ is a smooth irreducible divisor and $Z \subset H$, the complex log canonical threshold satisfies (see \cite{Mustata})
$$\text{lct}_Z(X,Y) \geq \text{lct}_Z(H,Y \cap H).$$
This inequality illustrates a well-known principle in complex analytic geometry: the singularities of the pair $( H, Y \cap H)$ are at least as bad as those of the pair $(X,Y)$. But this is not the case for the real log canonical threshold nor, a fortiori, for the real log canonical thresholds with signs.
    If one consider for instance $f(x,y,z,w) = x^{2}+y^{2}+z^{2}+w^{3}$ and $H= \{w=0\}$, then one has
    $$\text{rlct}_0(f) =1 \text{ but } \text{rlct}_0(f_{|H}) = \frac{3}{2}.$$
\end{rem}

\begin{rem}
    Assume that we are given a one-parameter family of analytic functions $f_t:(\mathbb{C}^d,0) \to (\mathbb{C},0)$ where $t$ belongs to an open set $U \subset \mathbb{C}$ such that the Milnor number $\mu(f_t)$ is constant for all $t\in U$. Then the log canonical threshold $\text{lct}_0(f_t)$ is also constant for all $t\in U$ (see \cite{varchenko2}), in other words every $\mu$-constant family is also $\text{lct}_0$-constant. The real log canonical threshold does not satisfy such a rigid property either. For example, consider the family $f_t(x,y,z) = x^{2}+y^{2}+z^{2}+tz$ which is such that for every $t \in \mathbb{R}$,
    $$\mu(f_t) = \dim_{\mathbb{C}} \mathbb{C}[x,y,z]/(\partial{f_t}/\partial x , \partial{f_t}/\partial y ,\partial{f_t}/\partial z) = 1.$$
On the other hand, one finds that $\text{rlct}_0(f_0) = \frac{3}{2}$ whereas $\text{rlct}_0(f_t) = 1$ for $t \neq 0$.
\end{rem}

\begin{prop}
    Assume that $f:(\mathbb{R}^d,0) \to (\mathbb{R},0)$ is a real analytic function germ, and let $m=\text{mult}_0(f)$ denote the algebraic multiplicity of $f$ at the origin. One has the inequality
    $$\text{rlct}_0(f) \leq \frac{d}{m}.$$
\end{prop}

\begin{proof}
    Let us consider a log resolution $\sigma : X' \to \mathbb{A}^d_{\mathbb{R}} $ of $f$. We can then decompose $\sigma$ as the composition of finitely many blow-ups along smooth centers $\sigma : X' =X_k \to \dots \dots X_1 \to X_0 = \mathbb{A}^{d}_{\mathbb{R}}$, and we may assume that the first blow-up is at the origin, in other words, that the morphism $X_1 \to X_0$ is given by the blow-up $\text{Bl}_0(\mathbb{A}^d_{\mathbb{R}}) \to \mathbb{A}^d_{\mathbb{R}}$. One can check that the exceptional divisor $E_1 \subset \sigma^{-1}(0)$ created by this blow-up has numerical data $(\nu_1,N_1)=(d,m)$ and that the real locus of 
    $$E_1\cap \sigma^{-1}(0)  = E_1 \simeq \mathbb{P}^{d-1}_{\mathbb{R}}$$
    is non-empty. Since 
    $$\text{rlct}_0(f) = \text{min} \{ \frac{\nu_i}{N_i} \mid i \in J, (E_i \cap \sigma^{-1}(0))(\mathbb{R}) \neq \emptyset  \}$$
    one obtains the desired inequality. Let us mention that one can also prove this inequality using the characterization of $\text{rlct}_0(f)$ in terms of jet spaces (cf. Proposition \ref{prop 1}), as is done in \cite{Mustata}. 
\end{proof}

\begin{rem}
    However, the two inequalities 
$$\text{rlct}_0^{+}(f) \leq \frac{d}{m} \text{ et } \text{rlct}_0^{-}(f) \leq \frac{d}{m}$$
do not always hold. Consider for example, for $k \geq 3$, the germ $f_k(x,y) = y^{8} + x^{4}y^{2} +x^{2k+1}$. One finds that (see Example \ref{ex 2} for the details)
$$\text{rlct}_0(f_k)=\text{rlct}^{+}_0(f_k)=\frac{5}{16} \leq \frac{1}{3} = \frac{d}{m}$$
while
$$\text{rlct}_0^{-}(f_k) = \frac{2k-1}{4k+2} > \frac{1}{3}.$$
\end{rem}

\begin{rem}
   Let us mention again that the real log canonical threshold is an invariant that takes into account the singularities of the real locus of an algebraic variety defined over $\mathbb{R}$. In some cases, it allows one to distinguish between two non-isomorphic real forms of a singular complex algebraic variety. Consider, for example, the varieties
 $$Y = \text{Spec}(\frac{\mathbb{R}[x,y,z]}{(x^{2}+y^{2}+z^{2})}) \subset \mathbb{A}^3_{\mathbb{R}} \text{ and } W = \text{Spec}(\frac{\mathbb{R}[x,y,z]}{(xy+z^{2})}) \subset \mathbb{A}^3_{\mathbb{R}}$$
 which satisfy
 $$Y_{\mathbb{C}} \simeq W_{\mathbb{C}}.$$
 On the other hand, one check that $\text{rlct}(\mathbb{A}^3_{\mathbb{R}},Y) = \frac{3}{2}$, whereas $\text{rlct}(\mathbb{A}^3_{\mathbb{R}},W) = 1$, therefore $Y$ and $W$ are not isomorphic over $\mathbb{R}$. Of course, there are many other ways to see this, such as considering the topological type of $Y(\mathbb{R})$ and $W(\mathbb{R})$, or computing the virtual Poincaré polynomial of these two varieties.
\end{rem}

\section{The Newton non-degenerate case}

\label{section 4}

In this section, we compute the real log canonical threshold and the log canonical thresholds with signs of  Nash functions germs $f : (\mathbb{R}^d,0) \to (\mathbb{R},0)$ which are Newton non-degenerate. Although the authors do not explicitly mention the log canonical threshold, results of this type already appear in \cite{goulw4}, where they are proved using a different approach and used to show that the weights are invariants of the blow-Nash equivalence for real quasi-homogeneous polynomials in three variables.

\begin{defi}
    Let $f:(\mathbb{R}^d,0) \to (\mathbb{R},0)$ be a Nash function germ, written in a system of local coordinates as
$$f(x)=\sum\limits_{i \in \mathbb{N}^d} a_i x^{i}$$
and denote by $\text{supp}(f) = \{ i \in \mathbb{N}^d \mid a_i \neq 0\}$ the support of $f$. We denote by $\Gamma(f)$ the Newton polyhedron of $f$, which is defined as the convex hull of the set
$$\text{supp}(f) + (\mathbb{R}_+)^{d}.$$ 
We also denote by $\Gamma_c(f)$ the union of the compact faces of $\Gamma(f)$, and we write $\tau \prec \Gamma_c(f)$ to indicate that $\tau$ is a compact face of $\Gamma(f)$. Finally, we say that $\tau$ is a facet if it is a compact codimension-one face of $\Gamma(f)$.
\end{defi}

\begin{defi}
     If $\tau \prec \Gamma_c(f)$, we denote by $f_{\tau}$ the quasi-homogeneous polynomial defined by
$$f_{\tau}(x)=\sum\limits_{i \in \tau} a_i x^{i}.$$
We say that $f$ is Newton non-degenerate if for every $\tau \prec \Gamma_c(f)$,
$$\{ x \in (\mathbb{R}^{*})^d \mid f_{\tau}(x) = \frac{\partial f_{\tau}}{\partial x_1}(x) = \dots = \frac{\partial f_{\tau}}{\partial x_d}(x) = 0 \} = \emptyset.$$
\end{defi}

\begin{defi}
    We say that $f$ is convenient if, for every $i \in \llbracket 1, d \rrbracket$, there exists an integer $p_i \in \mathbb{N}^*$ such that the monomial $x_i^{p_i}$ appears in the expansion of $f$.
\end{defi}

\begin{defi}
    For $k \in (\mathbb{R}_+)^d$, we set
$$\nu(k)=\langle k, (1,\dots,1) \rangle = k_1 +\dots +k_d.$$
We denote by $N$ the support function associated to $\Gamma(f)$, which is defined for $k \in (\mathbb{R}_+)^d$ by
$$N(k) = \text{inf} \{ \langle k ,x \rangle \mid x \in \Gamma(f) \} = \text{min} \{ \langle k ,x \rangle \mid x \in \Gamma(f)\}. $$
If $k \in (\mathbb{R}_+)^d$, we denote by $\tau(k)$ the face of $\Gamma_+(f)$ on which $N(k)$ is attained, that is,
$$\tau(k) = \{ x \in \Gamma(f) \mid  \langle k,x \rangle = N(k) \}$$
and $\tau(k)$ is then a compact face of $\Gamma(f)$ if and only if $k \in (\mathbb{R}^*_+)^d$.
Finally, if $\tau$ is a face of $\Gamma(f)$, we denote by $\sigma(\tau)$ the associated cone, defined by
$$\sigma(\tau) = \{ k \in (\mathbb{R}_+)^d \mid \tau = \tau(k) \}.$$
The union of these cones forms a fan $\Sigma_0$, which is dual to $\Gamma(f)$.
\end{defi}

\begin{thm}[\cite{varchenko}]
\label{thm 9}
  There exists a regular fan $\Sigma$ refining $\Sigma_0$, which induces a non-singular real toric variety $X(\Sigma)$ together with a proper morphism $\sigma : X(\Sigma) \to \mathbb{A}^d_{\mathbb{R}}$. The following conditions are then equivalent:
\begin{enumerate}
    \item The morphism $\sigma : X(\Sigma) \to \mathbb{R}^d$ is a log-resolution of the germ $f :(\mathbb{R}^d,0) \to (\mathbb{R},0)$.
    \item The germ $f : (\mathbb{R}^d,0) \to (\mathbb{R},0)$ is Newton non-degenerate.
\end{enumerate}
\end{thm}

\begin{defi}
   Let $\tau \prec \Gamma_c(f)$. We then define the real algebraic varieties $X_{\tau}, X_{\tau}^c$ and $X_{\tau}^{\pm}$ by the formulas
$$X_{\tau} = \{ x \in (\mathbb{R}^*)^d \mid f_{\tau}(x) \neq 0 \},$$

$$X_{\tau}^c = (\mathbb{R}^*)^d \setminus X_{\tau} = \{ x \in (\mathbb{R}^*)^d \mid f_{\tau}(x) = 0 \}$$
and
$$X_{\tau}^{\pm} = \{ x \in (\mathbb{R}^*)^d \mid f_{\tau}(x) = \pm 1 \}.$$
We also denote by $S_{\tau}$ the formal series defined by
$$S_{\tau} = \sum\limits_{ k \in \sigma(\tau) \cap (\mathbb{N}^*)^d} u^{-\nu(k)} T^{N(k)} \in \mathbb{Z}[u,u^{-1}][[T]].$$
\end{defi}

\begin{rem}
\label{rem 4}
   \begin{enumerate}
       \item For $\tau \prec \Gamma_c(f)$, $X_{\tau}$ is a Zariski open subset of $(\mathbb{R}^*)^d$, hence it is a real algebraic variety of dimension $d$. On the other hand, it may happen that $X_{\tau}^{c}$ is empty; however, if this is not the case, and under the assumption that $f$ is Newton non-degenerate, $X_{\tau}^{c}$ is then a non-singular real algebraic variety of dimension exactly $d-1$. By quasi-homogeneity, the varieties $X_{\tau}^{\pm}$ are also non-singular, and therefore, whenever nonempty, they are real algebraic varieties of dimension exactly $d-1$.
       \item Let $\tau \prec \Gamma_c(f)$ be a face and assume that $X_{\tau}^c \neq \emptyset$. Then, under the assumption that $f$ is Newton non-degenerate, $X_{\tau}^c$ is a smooth variety and $f_\tau$ changes sign in a neighborhood of $X_{\tau}^c$. By quasi-homogeneity, this implies that the subvarieties $X_{\tau}^{\pm}$ are nonempty.
       \item The explicit computation of the sums $S_{\tau}$ in the form of a rational function, as in the Denef--Loeser formula (cf. Theorem \ref{thm 2}), can be carried out by subdividing $\sigma(\tau)$ as a union of regular cones, which amounts precisely to computing the toric log-resolution $\sigma: X(\Sigma) \to \mathbb{R}^d$.
   \end{enumerate}
\end{rem}

The computation of zeta functions and invariants associated with these zeta functions in the non-degenerate case has, for instance, already been carried out in \cite{loes2}, \cite{guibert}, \cite{goulw4}, \cite{campesato}, etc.

\begin{thm}[\cite{campesato} Theorem 5.15]
\label{thm 6}
    Assume that the Nash function germ $f:(\mathbb{R}^d,0) \to (\mathbb{R},0)$ is Newton non-degenerate. Then one has
$$Z_{\beta,0}(f;T) =  \sum\limits_{ \tau \prec \Gamma_c(f)}\big(\beta(X_{\tau}) +(u-1)\beta(X_{\tau}^c)\frac{u^{-1}T}{1-u^{-1}T}\big) S_{\tau}.$$
For the zeta functions with signs, one has
$$Z_{\beta,0}^{\pm}(f;T) =  \sum\limits_{ \tau \prec \Gamma_c(f)}\big(\beta(X_{\tau}^{\pm}) +\beta(X_{\tau}^c)\frac{u^{-1}T}{1-u^{-1}T}\big) S_{\tau}.$$
\end{thm}

\begin{rem}
    In \cite{campesato}, this equality is proved at the motivic level, that is in the ring $\text{K}_0(AS)$, and the desired equality is then obtained by applying the virtual Poincaré polynomial. On the other hand, the formula is only proved for the naive zeta function $Z_{\beta,0}(f;T)$, but the argument can be easily adapted to prove the formula in the case of zeta functions with signs.
\end{rem}

We will always assume in the sequel of this section that the germ $f :(\mathbb{R}^d,0) \to (\mathbb{R},0)$ is Newton non-degenerate.

\begin{defi}
    We define $l_f,l_f^{\pm} \in \mathbb{Q}^*_+$ by the formulas
$$l_f = \text{inf} \{ \frac{\nu(k)}{N(k)} \mid \tau \prec \Gamma_c(f) \text{ and } k \in \sigma(\tau) \cap (\mathbb{N}^*)^d  \}  $$
and
$$l_f^{\pm} = \text{inf} \{ \frac{\nu(k)}{N(k)} \mid \tau \prec \Gamma_c(f),~X_{\tau}^{\pm} \neq \emptyset \text{ and } k \in \sigma(\tau) \cap (\mathbb{N}^*)^d\}.$$
\end{defi}

\begin{rem}
\label{rem 5}
    Since every face is an intersection of facets, one has
    
    $$ \{ \frac{\nu(k)}{N(k)} \mid \tau \prec \Gamma_c(f) \text{ and } k \in \sigma(\tau) \cap (\mathbb{N}^*)^d  \} = \{ \frac{\nu(k)}{N(k)} \mid \tau \text{ is a compact facet and } k \in \sigma(\tau) \cap (\mathbb{N}^*)^d  \}.$$
   Now, let us fix a compact facet $\tau$. The set $\sigma(\tau)$ is a one-dimensional cone so there exists a vector $a=(a_1, \dots,a_d) \in (\mathbb{N}^*)^d$ such that $\sigma(\tau) = \mathbb{R}_+a$, and we may also assume that $\gcd(a_1,\dots,a_d)=1$. Thus, every vector $k \in \sigma(\tau) \cap (\mathbb{N}^*)^d$ can be written $k = \lambda a $ for some $\lambda \in \mathbb{N}^*$ and therefore satisfies $\frac{\nu(k)}{N(k)} = \frac{\lambda \nu(a)}{\lambda N(a)} = \frac{\nu(a)}{N(a)}$. In particular
    $$\text{inf} \{ \frac{\nu(k)}{N(k)} \mid k \in \sigma(\tau) \cap (\mathbb{N}^*)^d  \} = \text{min} \{ \frac{\nu(k)}{N(k)} \mid k \in \sigma(\tau) \cap (\mathbb{N}^*)^d  \} = \frac{\nu(a)}{N(a)}.$$
    It follows that the infimum defining $l_f$ is in fact a minimum, and the same holds for $l_f^{\pm}$.
\end{rem}

\begin{thm}
\label{thm 10}
   Assume that $X_{\tau}^{c} = \emptyset$ for all $\tau \prec \Gamma_c(f)$. Then
$$\text{rlct}_0^{\pm}(f) = l_f^{\pm}.$$
Now assume that there exists at least one $\tau \prec \Gamma_c(f)$ such that $X_{\tau}^{c} \neq \emptyset$. Then
$$\text{rlct}_0^{\pm}(f) = \text{min}(1,l_f^{\pm}).$$
\end{thm}

\begin{proof}
We will use the following characterization of the real log-canonical threshold with signs (see Lemma \ref{lemme 1}):
$$\text{rlct}_0^{\pm}(f) = \text{sup } \{\alpha \in \mathbb{Q} \mid \text{ Every coefficient of } Z_{\beta,0}^{\pm}(f;u^\alpha T)u^{1-d} \text{ has degree } < 0\}.$$

Assume first that $X_{\tau}^{c} = \emptyset$ for all $\tau \prec \Gamma_c(f)$. Then Theorem \ref{thm 6} gives
$$Z_{\beta,0}^{\pm}(f;u^{\alpha}T)u^{1-d} =  \sum\limits_{ \tau \prec \Gamma_c(f)}\beta(X_{\tau}^{\pm})u^{1-d} S_{\tau}(u^{\alpha}T).$$
Moreover, when $X_{\tau}^{\pm} = \emptyset$, we have $\beta(X_{\tau}^{\pm}) = 0$, so that in the sum above we may restrict ourselves to faces $\tau \prec \Gamma_c(f)$ for which $X_{\tau}^{\pm} \neq \emptyset$.
Let us fix $\tau \prec \Gamma_c(f)$ such that $X_{\tau}^{\pm} \neq \emptyset$. By point (1) of Remark \ref{rem 4}, the term $\beta(X_{\tau}^{\pm})u^{1-d}$ is an element of $\bigcup_{k \geq 1}\mathbb{Z}[u^{\frac{1}{k}}, u^{-\frac{1}{k}}]$ of degree $0$ whose leading coefficient is strictly positive, and each coefficient of $S_{\tau}(u^{\alpha}T)$ is a finite sum of monomials with rational exponents. It follows that all coefficients of $Z_{\beta,0}^{\pm}(f;u^{\alpha}T)u^{1-d}$ have degree $< 0$ if and only if, for every $\tau \prec \Gamma_c(f)$ such that $X_{\tau}^{\pm} \neq \emptyset$, all coefficients of $S_{\tau}(u^{\alpha}T)$ have degree $< 0$. This holds if and only if, for every $\tau \prec \Gamma_c(f)$ such that $X_{\tau}^{\pm} \neq \emptyset$ and every $k \in \sigma(\tau) \cap (\mathbb{N}^*)^d$, the inequality $-\nu(k) +\alpha N(k) < 0$ holds. We thus obtain
$$\text{rlct}_0^{\pm}(f) = \text{inf} \left\{ \frac{\nu(k)}{N(k)} \mid \tau \prec \Gamma_c(f),~ X_{\tau}^{\pm} \neq \emptyset \text{ and } k \in \sigma(\tau) \cap (\mathbb{N}^*)^d \right\} = l_f^{\pm}.$$

Now assume that there exists $\tau \prec \Gamma_c(f)$ such that $X_{\tau}^{c} \neq \emptyset$. We can decompose $Z_{\beta,0}^{\pm}(f;T)u^{1-d}$ as
$$Z_{\beta,0}^{\pm}(f,u^{\alpha}T)u^{1-d} = \sum\limits_{ \tau \prec \Gamma_c(f)}\beta(X_{\tau}^{\pm})u^{1-d} S_{\tau}(u^{\alpha}T) + \sum\limits_{ \tau \prec \Gamma_c(f)}\beta(X_{\tau}^{c})u^{1-d} \frac{u^{\alpha-1}T}{1-u^{\alpha-1}T} S_{\tau}(u^{\alpha}T)$$
where in first sum we may restrict to faces $\tau \prec \Gamma_c(f)$ for which $X_{\tau}^{\pm} \neq \emptyset$ and in the second sum we may restrict to faces $\tau \prec \Gamma_c(f)$ for which $X_{\tau}^{c}$ is non-empty, and the latter contains at least one term by assumption.

Let us fix $\tau \prec \Gamma_c(f)$ such that $X_{\tau}^c \neq \emptyset$. By point (1) of Remark \ref{rem 4}, the term $\beta(X_{\tau}^{c})u^{1-d}$ is an element of $\bigcup_{k \geq 1}\mathbb{Z}[u^{\frac{1}{k}}, u^{-\frac{1}{k}}]$ of degree $0$ whose leading coefficient is strictly positive, and each coefficient of $\frac{u^{\alpha-1}T}{1-u^{\alpha-1}T}S_{\tau}(u^{\alpha}T)$ is a sum of monomials with rational exponents. Therefore all coefficients of $Z_{\beta,0}^{\pm}(f;u^{\alpha}T)u^{1-d}$ have degree $< 0$ if and only if, for every $\tau \prec \Gamma_c(f)$ such that $X_{\tau}^{\pm} \neq \emptyset$, all coefficients of $S_{\tau}(u^{\alpha}T)$ have degree $0$ and for every $\tau \prec \Gamma_c(f)$ such that $X_{\tau}^c \neq \emptyset$, all coefficients of $\frac{u^{\alpha-1}T}{1-u^{\alpha-1}T}S_{\tau}(u^{\alpha}T)$ have degree $< 0$. This holds if and only if, for every $\tau \prec \Gamma_c(f)$ such that $X_{\tau}^{\pm} \neq \emptyset$ and for every $k \in \sigma(\tau) \cap (\mathbb{N}^*)^d$, one has $-\nu(k) +\alpha N(k) < 0$ and for every $\tau \prec \Gamma_c(f)$ such that $X_{\tau}^{c} \neq \emptyset$ and for every $k \in \sigma(\tau) \cap (\mathbb{N}^{*})^{d}$ and $ m \geq 1$, one has  $-m +\alpha m -\nu(k) +\alpha N(k) < 0$. By point (2) of Remark \ref{rem 4}, this holds if and only if for every $\tau \prec \Gamma_c(f)$ such that $X_{\tau}^{\pm} \neq \emptyset$, for every $k \in \sigma(\tau) \cap (\mathbb{N}^{*})^d$ one has $-\nu(k)+\alpha N(k) <0 $ and for every $m \geq 1$ one has $-m + \alpha m < 0$. We thus obtain
$$\text{rlct}_0^{\pm}(f) = \text{min}(1, l_f^{\pm}).$$
\end{proof}

\begin{rem}
    Consider the toric log-resolution $\sigma: X(\Sigma) \to \mathbb{A}^d_{\mathbb{R}}$ of Theorem \ref{thm 9}. Let us fix a compact face $\tau \prec \Gamma_c(f)$ and consider $X(\sigma(\tau)) \subset X(\Sigma)$ the associated affine toric variety. Then the results of \cite{varchenko} also show that the equation of the strict transform in the open set $X(\sigma(\tau))$ is exactly given by $f_{\tau} = 0$. In other words, all the sets $X_{\tau}^{c}$ are empty if and only if the strict transform of the germ $f :(\mathbb{R}^d,0) \to (\mathbb{R},0)$ has no real points.

Note that since $f$ is non-degenerate, $f$ is reduced, so that $\frac{\nu_i}{N_i}= 1$ whenever $E_i$ is a component of the strict transform, which is consistent with the expression of the real log-canonical threshold given above.
\end{rem}

\begin{thm}
Assume that $X_{\tau}^{c} = \emptyset$ for all $\tau \prec \Gamma_c(f)$. Then
$$\text{rlct}_0(f) = l_f.$$
Now assume that there exists at least one $\tau \prec \Gamma_c(f)$ such that $X_{\tau}^{c} \neq \emptyset$. Then
$$\text{rlct}_0(f) = \text{min}(1,l_f).$$
\end{thm}

\begin{proof}
   The proof is the same as that of the previous theorem, using the characterization
$$\text{rlct}_0(f) = \text{sup } \{\alpha \in \mathbb{Q} \mid \text{ Every coefficient of } Z_{\beta,0}(f;u^\alpha T)u^{-d} \text{ has degree } < 0\}$$
and the formula for naive zeta functions from Theorem \ref{thm 6}.
\end{proof}

One can also read off the numbers $l_f,l_f^{\pm}$ quite easily from the Newton polyhedron of $f$.
Fix a compact facet $\tau$, so that $\sigma(\tau)$ is a one-dimensional cone, and denote by $H_{\tau}$ the supporting hyperplane of $\tau$. We then define the quantities
$$l_{f,\tau} = \text{inf} \left\{ \frac{\nu(k)}{N(k)} \mid k \in \sigma(\tau) \cap (\mathbb{N}^*)^d \right\}$$
and
$$d_{f,\tau} = \text{inf} \left\{ c >0 \mid c(1,\dots,1) \in H_{\tau} \right\}.$$
Moreover, as shown in Remark \ref{rem 5}, these infima are in fact minima.
\begin{lem}
    For every compact facet $\tau$, one has $d_{f,\tau} = \frac{1}{l_{f,\tau}}.$
\end{lem}

\begin{proof}
    Let us fix a compact facet $\tau$. For every $k \in \sigma(\tau)$, an equation of the hyperplane $H_{\tau}$ is given by
    $$H_\tau = \{ x \in \mathbb{R}^d \mid \langle k ,x \rangle = N(k) \}.$$
Since $d_{f,\tau}(1,\dots,1) \in H_{\tau}$, one deduces in particular that $\langle k , d_{f,\tau}(1,\dots,1) \rangle = N(k)$, i.e., that $d_{f,\tau}\nu(k) = N(k)$. On the other hand, $l_{f,\tau}$ is attained, in other words there exists some $k \in \sigma(\tau)$ such that
$$l_{f,\tau} = \frac{\nu(k)}{N(k)} = \frac{1}{d_{f,\tau}}.$$
\end{proof}

\begin{prop}
\label{prop 6}
One has 
$$\frac{1}{l_f^{\pm}} = \text{max} \left\{ d_{f,\tau} \mid \tau \text{ is a compact facet and } X_{\tau}^{\pm} \neq \emptyset \right\}$$
and
$$\frac{1}{l_f} = \text{max} \left\{ d_{f,\tau} \mid \tau \text{ is a  compact facet} \right\} = \text{min} \left\{ c >0 \mid c(1,\dots,1) \in \Gamma_c(f) \right\}$$
If $f$ is convenient, one has in particular that
$$\frac{1}{l_f} = \text{min} \left\{ c >0 \mid c(1,\dots,1) \in \Gamma(f) \right\}.$$
\end{prop}

\begin{proof}
First, note that every compact face of $\Gamma(f)$ is an intersection of facets, so that
$$l_f^{\pm} = \text{min} \left\{ \frac{\nu(k)}{N(k)} \mid \tau \text{ is a facet },~X_{\tau}^{\pm} \neq \emptyset \text{ and } k \in \sigma(\tau) \cap (\mathbb{N}^*)^d\right\}$$
and
$$l_f = \text{min} \left\{ \frac{\nu(k)}{N(k)} \mid \tau \text{ is a facet and } k \in \sigma(\tau) \cap (\mathbb{N}^*)^d  \right\}.$$
By the previous lemma, it follows that
$$\frac{1}{l_f^{\pm}} = \frac{1}{\text{min} \left\{ l_{f,\tau} \mid \tau \text{ is a facet and }X_{\tau}^{\pm} \neq \emptyset \right\}} = \text{max} \left\{ d_{f,\tau} \mid \tau \text{ is a facet and } X_{\tau}^{\pm} \neq \emptyset \right\}$$
and similarly that
$$\frac{1}{l_f} = \frac{1}{\text{min} \left\{ l_{f,\tau} \mid \tau \text{ is a facet} \right\}} =  \text{max} \left\{ d_{f,\tau} \mid \tau \text{ is a facet} \right\}.$$
By convexity of $\Gamma_c(f)$ which lies above all its supporting hyperplanes, we now have 
$$\text{min} \left\{ c >0 \mid c(1,\dots,1) \in \Gamma_c(f) \right\} =  \text{max} \left\{ d_{f,\tau} \mid \tau \text{ is a facet} \right\} = \frac{1}{l_f}. $$
Assume now that $f$ is convenient. In particular, every noncompact face of $\Gamma(f)$ is contained in a coordinate hyperplane. It follows that 
$$\text{min} \left\{ c >0 \mid c(1,\dots,1) \in \Gamma(f) \right\} = \text{min} \left\{ c >0 \mid c(1,\dots,1) \in \Gamma_c(f) \right\}.$$
\end{proof}

\begin{ex} 
1) Consider $f(x,y,z) = x^{2}+y^{2}+z^{2k} +xyz$ with $k \geq 1$. Then the Newton polyhedron of $f$ has a single face $\tau$ passing through the points $(2,0,0)$, $(0,2,0)$ and $(0,0,2k)$, and which is therefore supported by the hyperplane of equation $kx+ky+z=2k$. It follows that
$$l_f = \frac{2k+1}{2k}$$
and since
$$X_{\tau}^c = \{(x,y,z) \in (\mathbb{R}^*)^{3} \mid x^{2} + y^{2} +z^{2k} = 0 \} = \emptyset$$
one has
$$\text{rlct}_0(f) = l_f = \frac{2k+1}{2k} \text{ while } \text{lct}_0(f) = 1.$$

\noindent
2) More generally, assume that $f$ is Newton non-degenerate and convenient, and that $f$ is semi-quasihomogeneous. In other words, one can write $f= h+g $ where $h$ is a quasihomogeneous polynomial associated with a weight vector $(\omega ;\omega_1,\dots,\omega_d)$ and $g$ is a finite linear combination of monomials of degree $> \omega$ with respect to the weight vector $(\omega_1,\dots,\omega_d)$. This means precisely that $\Gamma(f)$ has a unique compact face $\tau$ supported by the hyperplane $H_{\tau}$ with equation
$$\omega_1x_1+\dots+\omega_dx_d = \omega$$
and that
$f_{\tau} = h$. It follows that if
$$\{ x \in (\mathbb{R}^*)^d \mid h(x) = 0 \} = \emptyset,$$
then $\text{rlct}_0(f) = \frac{\omega_1+\dots+\omega_d}{\omega}$, whereas otherwise one has $\text{rlct}_0(f) = \text{min}(1,\frac{\omega_1+\dots+\omega_d}{\omega})$.

\end{ex}

\begin{ex}
   Take $f(x,y) = x^9+y^4+x^4y^2$, which is Newton non-degenerate and convenient. The Newton polygon of $f$ has two compact facets $\tau_1,\tau_2$ and is depicted below.
    \begin{figure}[H]
        \centering
        \includegraphics[scale=0.35]{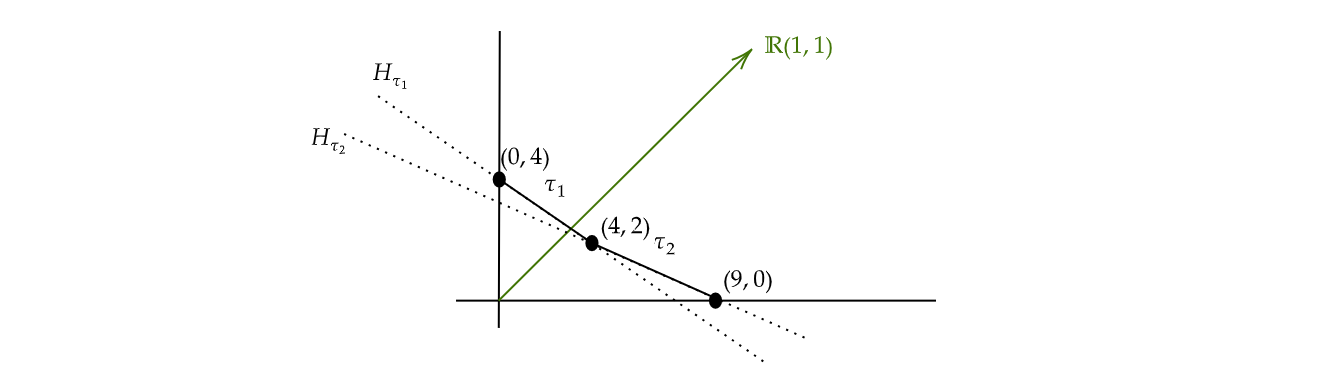}
    \end{figure}
   The supporting lines $H_{\tau_1}$ and $H_{\tau_2}$ have equations $2y+x=8$ and $5y+2x=18$, respectively. One sees that the line $\mathbb{R}(1,1)$ intersects $\Gamma_c(f)$ on the facet $\tau_1$ at the point $(\frac{8}{3},\frac{8}{3})$, so that $\text{rlct}_0(f) = \text{rlct}_0^+(f) = \frac{3}{8}$. On the other hand, one finds that
$$X_{\tau_1}^{-} = \{ (x,y) \in (\mathbb{R}^*)^2 \mid x^9+x^4y^2= -1 \}  =\emptyset,$$
whereas $X_{\tau_2}^- \neq \emptyset$. The line $\mathbb{R}(1,1)$ intersects $H_{\tau_2}$ at the point $(\frac{18}{7},\frac{18}{7})$, and it follows that $\text{rlct}_0^-(f) = \frac{7}{18} > \frac{3}{8} =\text{rlct}_0^+(f)$. One may also check that the graph of the associated embedded toric resolution is as shown below.

    \begin{figure}[H]
        \centering
        \includegraphics[scale=0.35]{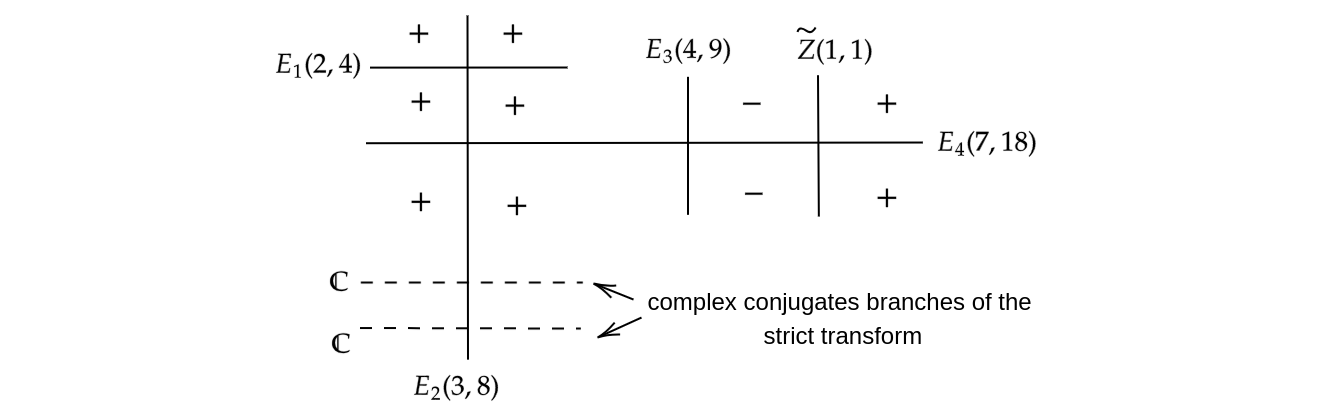}
    \end{figure}

\end{ex}

\begin{ex}
\label{ex 2}
Consider, for $k \geq 3$, the germ $f_k(x,y) = y^{8} + x^{4}y^{2} +x^{2k+1}$, which is Newton non-degenerate and convenient. The Newton polygon of $f_k$ has two compact facets $\tau_1,\tau_2$ and is depicted below.

 \begin{figure}[H]
        \centering
        \includegraphics[scale=0.35]{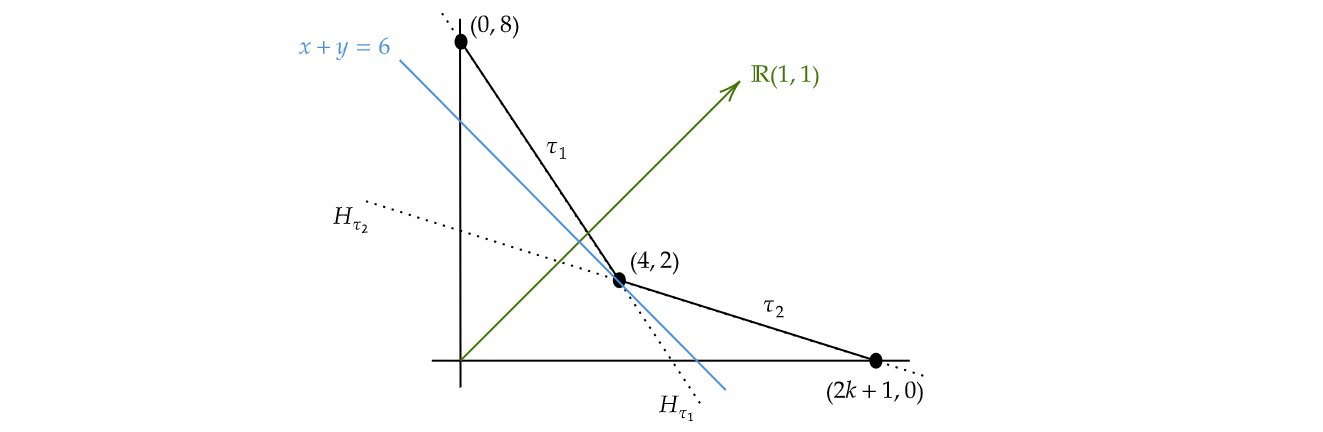}
    \end{figure}

The supporting lines $H_{\tau_1}$ and $H_{\tau_2}$ have equations $2y+3x=16$ and $(2k-3)y+2x=4k+2$, respectively. The line $\mathbb{R}(1,1)$ intersects the Newton polygon of $f_k$ along the face $\tau_1$ at the point $(\frac{16}{5},\frac{16}{5})$. With the notation and results of Section \ref{section 2}, one has
$$X_{\tau_1}^{+} = \{ (x,y) \in (\mathbb{R}^*)^{2} \mid y^{8}+x^{4}y^{2} = 1 \} \neq \emptyset$$
and therefore 
$$\text{rlct}_0(f_k)=\text{rlct}^{+}_0(f_k)=\frac{5}{16} \leq \frac{1}{3} = \frac{d}{m}.$$
On the other hand, one finds that 
$$X_{\tau_1}^{-} = \{(x,y) \in (\mathbb{R}^{*})^{2} \mid y^{8}+x^{4}y^{2} = -1 \} = \emptyset$$
and that the line $\mathbb{R}(1,1)$ intersects $H_{\tau_2}$ at the point $(\frac{4k+2}{2k-1}, \frac{4k+2}{2k-1})$ where
$$X_{\tau_2}^{-} = \{(x,y) \in (\mathbb{R}^*)^2 \mid x^4y^2+x^{2k+1} = -1 \} \neq \emptyset.$$
Thus
$$\text{rlct}_0^{-}(f_k) = \frac{2k-1}{4k+2} > \frac{1}{3}.$$

\end{ex}

\section{Plane curve singularities}

\label{section 5}
In this section we study real log canonical thresholds and log canonical thresholds with signs in the case of plane curves. In particular, we compute the pair $\{ \text{rlct}_0^{+}(f), \text{rlct}_0^{-}(f) \}$ in terms of characteristic Puiseux exponents in the irreducible case.

\begin{prop}
    Let $f :(\mathbb{R}^{2},0) \to (\mathbb{R},0)$ be a real analytic function germ. Then one has 
    $$\text{rlct}_0(f) = \text{lct}_0(f).$$
    
\end{prop}

\begin{proof}
    Let $\sigma : X \to \mathbb{A}^{2}_{\mathbb{R}}$ be a log resolution of $f$. We have seen that
$$\text{rlct}_0(f) = \text{min} \{ \frac{\nu_i}{N_i} \mid i \in J_{\mathbb{R}},~ 0\in \sigma(E_i)\} \text{ and } \text{lct}_0(f) = \text{min} \{ \frac{\nu_i}{N_i} \mid i \in J, \ 0 \in \sigma(E_i) \} $$
    If $\text{lct}_0(f)$ is realized by an exceptional curve $E_i \simeq \mathbb{P}^{1}_{\mathbb{R}}$, then $E_i(\mathbb{R}) \neq \emptyset$, and hence $E_i$ also realizes $\text{rlct}_0(f)$.
    Assume now that $\text{lct}_0(f)$ is realized by a component of the strict transform of $f$ and that for every exceptional curve $E_i$ one has $\frac{\nu_i}{N_i} > \text{lct}_0(f)$. We write $f = f_1^{M_1}\dots f_r^{M_r}$ the decomposition of $f$ into irreducible factors, and we therefore assume that $\text{lct}_0(f)=\frac{1}{M_j}$ for some component $f_j$ such that $f_j(0)=0$. For $k \in \llbracket 1, r \rrbracket$ we denote by $\mu_k = \text{mult}_0(f_k)$ the algebraic multiplicity of $f_k$ at the origin and we consider the first exceptional curve $E_1$, obtained by blowing up the origin. Then we have $\nu_1=2$ and $N_1 = \sum\limits_{k=1}^r \mu_k M_k$, and by assumption $\frac{2}{N_1} > \frac{1}{M_j}$. This forces $\mu_j =1$, so that the origin is a smooth point of $\{f_j=0\}$, and hence the real locus of $f_j$ is not reduced to the origin. Since $\sigma$ is an isomorphism outside the origin, this implies that the real locus of the strict transform of $f_j$ is nonempty, and hence that this component also realizes $\text{rlct}_0(f)$.  
\end{proof}

\begin{assumption}
\label{assumption 1}
    From now on, we assume that $f:(\mathbb{R}^{2},0) \to (\mathbb{R},0)$ is analytically irreducible over $\mathbb{R}$ and we denote $\sigma : X \to \mathbb{A}^2_{\mathbb{R}}$ the canonical embedded resolution of $f$. We also assume that the real locus of $f$ is not reduced to the origin, since otherwise one has
$$\text{rlct}_0^+(f)  = + \infty \text{ or } \text{rlct}_0^-(f)  = + \infty.$$
\end{assumption}

These two assumptions implies in particular that $f$ is also analytically irreducible in $\mathbb{C}\{x,y\}$.

\begin{defi}
   The total dual graph of the resolution is the graph $G$ whose vertices are the components $E_i$ for $i \in J$, and two vertices $i,j$ are connected by an edge if and only if $E_i \cap E_j \cap \sigma^{-1}(0) \neq \emptyset$. We represent by a point the vertices corresponding to exceptional curves and by a circle the unique vertex corresponding to the analytically irreducible component of the strict transform. 
\end{defi}

Under the assumptions made above, the total dual graph considered here coincides with the real total dual graph introduced in \cite{jaudon2026polesrealmotiviczeta}.

\begin{defi}
   We write the canonical embedded resolution of $f$ as the composition of a finite number of blow-ups
$\sigma : X =X_k \to \dots \dots X_1 \to X_0 = \mathbb{A}^{2}_{\mathbb{R}}$,
and for $i \in \llbracket 0, k-1 \rrbracket$ we denote by $f_i$ the strict transform of $f$ in $X_i$. Since $f$ is analytically irreducible, the zero locus of $f_i$ intersects the exceptional locus of $\sigma_i : X_i \to X_{i-1} \to \dots X_0$ at a unique point $x_i$ (that is, the next center of blow-up if $i < k$). We denote by $\mu_i = \text{mult}_{x_i}(f_i)$ and the multiplicity sequence of $f$ is defined to be $(\mu_0,\mu_1,\dots,\mu_{k-1})$.

\end{defi}

\begin{lem}
\label{lemme 2}
 For $i \in \llbracket 1, k \rrbracket$, let $E_i$ denote the exceptional divisor created by the blow-up $X_i \to X_{i-1}$. Then there exists a lower triangular matrix $M \in GL_k(\mathbb{Z})$ in which all the diagonal coefficients are 1 such that
$$(N_1,\dots,N_k)^{T} = M (\mu_0,\dots,\mu_{k-1})^{T}.$$

\end{lem}

\begin{proof}
  We argue by strong induction on $i\in \llbracket 1, k \rrbracket$, first noting that $N_1= \mu_0$. Now let $i \in \llbracket 2, k \rrbracket$. There are two possible cases. Either $E_i$ is created by blowing up an intersection point $E_{i-1} \cap E_j$ with $j < i$, in which case $N_i = N_{i-1} +N_j + \mu_{i-1}$. Or $E_i$ is created by blowing up a point lying only on $E_{i-1}$, in which case we have $N_i = N_{i-1} + \mu_{i-1}$. In both cases, it follows by induction that there exist integers $\lambda_0,\dots,\lambda_{i-2}$ such that
  $$N_i = \mu_{i-1} +\sum\limits_{l=0}^{i-2} \lambda_l \mu_l.$$

\end{proof}

\begin{notation}
   Consider a real Puiseux parametrization of $f$ of the form $t \to (t^m,\sum\limits_{i \geq n } \alpha_i t^i)$ with $m \leq n$. We also assume that it is a good parametrization, that is, there does not exist an integer $k \geq 2$ such that $k$ divides $m$ and $k$ divides $i$ for every $i \geq n$ such that $\alpha_i \neq 0$. We then denote by $(m;n_1,\dots,n_g)$ the associated characteristic Puiseux exponents, which are defined recursively as follows. We set
$$n_1 = \text{min} \{ i \mid \alpha_i \neq 0 \text{ and } m \nmid i\} $$
and for $ j \geq 2$ we set
$$n_j = \text{min} \{ i \mid \alpha_i \neq 0 \text{ and } \gcd(m,n_1,\dots,n_{j-1}) > \gcd(m,n_1,\dots,n_{j-1},i) \}.$$
Note that we then have $m = \text{mult}_0(f)$ and that $\gcd(m,n_1,\dots,n_g)=1$. For $1 \leq i \leq g$, we also denote by $e_i =\gcd(m,n_1,\dots,n_i)$, and we adopt the convention $e_0 = n_0 = m$.

\end{notation}

It is known (e.g. \cite{brieskorn}) that the total dual graph $G$ has the shape of a hanging mobile suspended from the component of the strict transform, with exactly $g$ rupture vertices (that is, vertices of valency 3) corresponding to exceptional divisors $E_{i_1},\dots,E_{i_g}$, as shown below.

    \begin{figure}[H]
        \centering
        \includegraphics[scale=0.35]{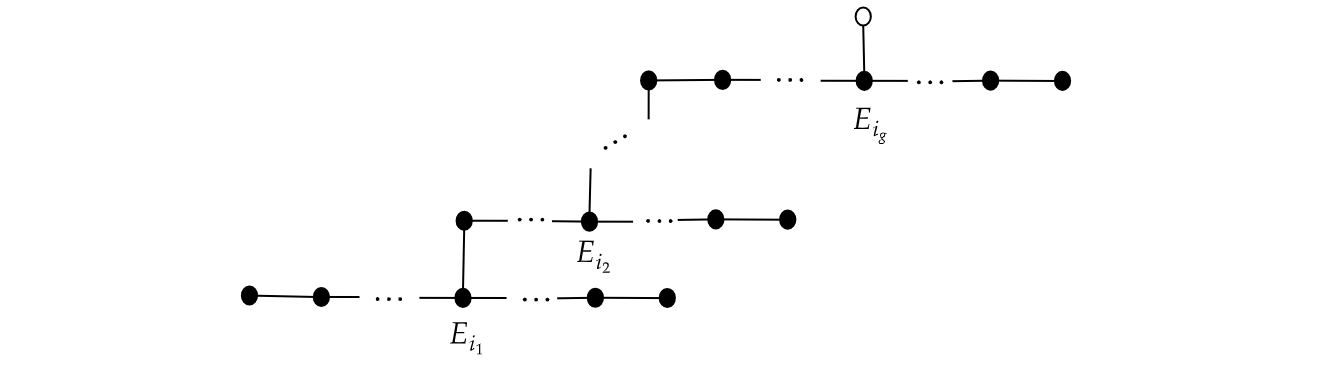}
    \end{figure}
By Theorem \ref{thm 3}, we know that $\text{rlct}^{\pm}_0(f)$ is a pole of $Z_{mot,0}^{\pm}(f;T)$, and by Corollary 3.17 of \cite{jaudon2026polesrealmotiviczeta} we know that for curves the poles of $Z_{mot,0}^{\pm}(f;T)$ are induced either by the exceptional divisors corresponding to rupture vertices or by components of the strict transform.

Moreover, it is also known (see \cite{strauss} or \cite{veys1}) that we have the following inequalities
$$\frac{\nu_{i_1}}{N_{i_1}} \leq \frac{\nu_{i_2}}{N_{i_2}} \leq \dots \leq \frac{\nu_{i_g}}{N_{i_g}} \leq 1.$$

A classical result of Igusa allows one to express the numerical data associated with the first rupture vertex in terms of the first two characteristic exponents via the formula
$$\text{lct}_0(f) = \frac{\nu_{i_1}}{N_{i_1}} = \frac{1}{m}+\frac{1}{n_1}.$$

More generally, there are formulas allowing one to compute the quotients $\frac{\nu}{N}$ for each rupture vertex in terms of the characteristic exponents.

\begin{thm}[\cite{Wall_2004} Theorem 8.5.2]
   We define integers $\overline{n}_0,\dots,\overline{n}_g$ recursively by setting $\overline{n}_0 = m$, $\overline{n}_1 = n_1$, and for $j \geq 2$,
$$\overline{n}_j = n_j - n_{j-1} + \frac{e_{j-2}}{e_{j-1}}\overline{n}_{j-1}.$$

Then for all $1 \leq j \leq g$, one has
$$\frac{\nu_{i_j}}{N_{i_j}} = \frac{n_j+ m}{e_{j-1} \overline{n}_j}.$$

\end{thm}

Our goal is now to compute the pair $\{\text{rlct}_0^{-}(f),\text{rlct}_0^{+}(f)\}$ in terms of the characteristic Puiseux exponents $(m;n_1,\dots,n_g)$.

To fix ideas, assume (up to replacing $f$ by $-f$) that $E_{i_1}(\mathbb{R})$ intersects $\overline{P(f)}$, and that we therefore have
$$\text{rlct}_0^{+}(f)=\text{rlct}_0(f)=\text{lct}_0(f) = \frac{\nu_{i_1}}{N_{i_1}}=\frac{1}{m} + \frac{1}{n_1}.$$

We already know that $E_{i_g}(\mathbb{R}) \cap \overline{N(f)} \neq \emptyset$, since $E_{i_g}$ intersects the strict transform of $f$, in a neighborhood of which $f\circ \sigma$ changes sign  because $f$ is assumed to be analytically irreducible. We can therefore define
$$h =\text{min} \{ 1 \leq j \leq g \mid E_{i_j}(\mathbb{R}) \cap \overline{N(f)} \neq \emptyset \}, $$
and from the discussion above we then have
$$\text{rlct}_0^{-}(f) = \frac{\nu_{i_h}}{N_{i_h}}=\frac{n_h+ m}{e_{h-1} \overline{n}_h}.$$

In the case of analytically irreducible curves, the data of the characteristic Puiseux exponents is an invariant from which many other well-known invariants can be expressed, such as the multiplicity sequence or the resolution graph. Our aim is now to characterize the integer h in terms of the Puiseux characteristic exponents. The following theorem will be useful in what follows, as it establishes a link between the Puiseux characteristic exponents and the sequence of multiplicities, and hence, by Lemma \ref{lemme 2}, between the Puiseux characteristic exponents and the multiplicities of the divisors appearing in the canonical embedded resolution.

\begin{thm}[\cite{brieskorn} section 8.4, Theorem 12]
\label{thm 7}
    The multiplicity sequence can be obtained from the characteristic exponents by performing $g$ Euclidean algorithms as follows. We set $s_1 = n_1$ and $s_j = n_j - n_{j-1}$ for $j \in \llbracket 2, g \rrbracket$. We then define numbers $(l_{j,i})$ and $(k_{j,i})$ with $1 \leq j \leq g$ and $1 \leq i \leq t_j$ by setting $l_{1,1} = m$ and then
$$s_j = k_{j,1}l_{j,1} + l_{j,2}  \text{ where } 0 \leq l_{j,2} < l_{j,1}$$
$$l_{j,1} = k_{j,2}l_{j,2} + l_{j,3} \text{ where } 0 \leq l_{j,3} < l_{j,2} $$
$$\vdots$$
$$l_{j,t_j-1} = k_{j,t_j}l_{j,t_j}$$
with $l_{j,1} = l_{j-1,t_{j-1}}$ for $j \in \llbracket 2, g \rrbracket$. The multiplicity sequence of $f$ is then given by
$$(l_{1,1},\dots,l_{1,1},\dots,l_{j,i},\dots,l_{j,i},\dots,1,\dots,1)$$
where each $l_{j,i}$ is repeated $k_{j,i}$ times.

\end{thm}

\begin{rem}
\label{rem 2}
   Note that Theorem \ref{thm 7} of \cite{brieskorn} is in fact more precise, as it also explains how this sequence of $g$ Euclidean algorithms allows one to reconstruct the dual resolution graph together with the order in which the exceptional divisors appear. For instance, let $\mathcal{P}_j \subset G$ denote the chain containing the $j$-th rupture vertex, as illustrated below.

\begin{figure}[H]
\centering
\includegraphics[scale=0.35]{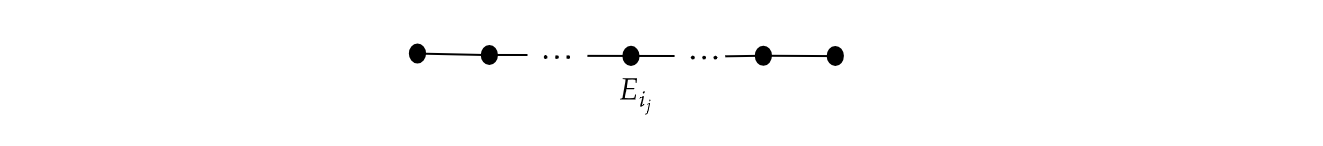}
\end{figure}

Then this chain is obtained by blowing up $k_{j,i}$ times, for all $i \in \llbracket 1, t_j \rrbracket$, the singular point where the strict transform has multiplicity $l_{j,i}$. In particular, the chain $\mathcal{P}_j$ contains $\sum\limits_{i=1}^{t_j} k_{j,i}$ vertices.

\end{rem}

\begin{lem}
    \begin{enumerate}
       \item Let $E_i$ be an exceptional curve, and assume that $E_i$ intersects other components $E_1,\dots,E_k$. Then one has 
$$-N_iE_i^2 = \sum\limits_{j=1}^k N_j .$$

\item Let $j\in \llbracket 1, g-1 \rrbracket$, and assume that $N_{i_j}$ is even and that for every component $E_{i}'$ intersecting $E_{i_j}$, the multiplicity $N_{i}'$ is even. Then for all $k \in \mathcal{P}_j$, the multiplicity $N_k$ is even.

    \end{enumerate}
\end{lem}

\begin{proof}
   $(1)$ Let us first note that, by Assumption \ref{assumption 1}, all intersection points of $E_i$ with the other components $E_1,\dots,E_k$ are real. By intersecting the divisor $E_i$ with the principal divisor $\text{div}(f \circ \sigma)$ in $\text{Pic}(X)$, we obtain
   $$0 = (E_i \cdot\sum\limits_{j \in J}N_j E_j) = N_iE_i^{2}  + \sum\limits_{j=1}^k N_j$$
   which yields the desired equality.
\newline
$(2)$ This point is proved by induction on the length of the chain $P_j$. The situation is as follows

\begin{figure}[H]
\centering
\includegraphics[scale=0.35]{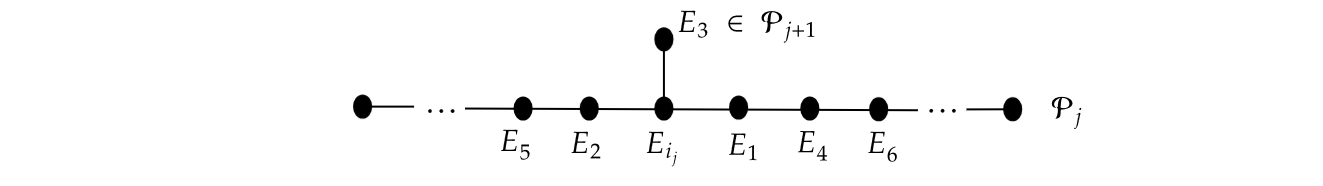}
\end{figure}

where $N_{i_j}, N_1,N_2,N_3$ are even. By applying point $(1)$ to the divisor $E_1$ one has
$$-N_1E_1^2 = N_{i_j} + N_4,$$
and it follows that $N_4$ is even. Similarly, by applying point $(1)$ to the divisor $E_2$, one finds that $N_5$ is even. Once again, by applying point $(1)$ to the divisor $E_4$, one finds that $N_6$ is even and more generally the desired result follows by induction on the length of the chain $\mathcal{P}_j$.  
\end{proof}

As we did before, let us assume that $E_{i_1}(\mathbb{R})$ intersects $\overline{P(f)}$ and consider the integer
    $$h = \text{min} \{ 1 \leq j \leq g \mid E_{i_j}(\mathbb{R}) \cap \overline{N(f)} \neq \emptyset \}.$$

\begin{prop}
\label{prop 5}
   The integer $h$ can be characterized by the following two equivalent conditions.
\begin{enumerate}
\item For all $j \in \underset{1 \leq i \leq h-1}{\bigcup}\mathcal{P}_i$, $N_j$ is even, and there exists $j \in \mathcal{P}_h$ such that $N_j$ is odd.
\item For all $j \in \llbracket1, h-1 \rrbracket$ and for all $i \in \llbracket 1, t_j \rrbracket$, $l_{j,i}$ is even, and there exists $i \in \llbracket 1, t_h \rrbracket$ such that $l_{h,i}$ is odd.
\end{enumerate}

\end{prop}

\begin{proof}
 Consider the chain $\mathcal{P}_1$ containing the first rupture vertex, where we recall that we have assumed that $E_{i_1}(\mathbb{R})\cap \overline{P(f)} \neq \emptyset$. Therefore $E_{i_1}(\mathbb{R}) \cap \overline{N(f)} = \emptyset$ if and only if $N_{i_1}$ is even and $N_{i}'$ is even for every component $E_{i}'$ intersecting $E_{i_1}$. By the previous lemma, this is also equivalent to the fact that for all $j \in \mathcal{P}_1$, $N_j$ is even. Applying the same argument to the subsequent chains we obtain by induction that the integer $h$ is characterized by condition $(1)$.
\newline
$(1) \iff (2)$ follows from Lemma \ref{lemme 2} and Remark \ref{rem 2}.

\end{proof}

\begin{thm}
\label{thm 8}
   If $m$ is odd, then one has
$$\text{rlct}_0^{-}(f)=\text{rlct}_0^{+}(f)= \frac{\nu_{i_1}}{N_{i_1}}=\frac{1}{m} + \frac{1}{n_1}.$$

If $m$ is even, then
$$\{ \text{rlct}_0^{+}(f),\text{rlct}_0^{-}(f) \} = \{ \frac{\nu_{i_1}}{N_{i_1}} , \frac{\nu_{i_h}}{N_{i_h}} \} = \{ \frac{1}{m} + \frac{1}{n_1} , \frac{n_h+ m}{e_{h-1} \overline{n}_h} \}$$ where
$$h = \text{min} \{ 1 \leq j \leq g \mid n_j \text{ is odd} \}.$$

\end{thm}

\begin{proof}
    Without loss of generality, after replacing $f$ with $-f$ if necessary (which exchange the values of $\text{rlct}_0^{+}(f)$ and $\text{rlct}_0^{-}(f)$), we may assume as before that $E_{i_1}(\mathbb{R}) \cap \overline{P(f)} \neq \emptyset$ and recall that $h$ is defined by
$$h =\text{min} \{ 1 \leq j \leq g \mid E_{i_j}(\mathbb{R}) \cap \overline{N(f)} \neq \emptyset \}.$$
    We have seen that 
    $$(\text{rlct}_0^{+}(f),\text{rlct}_0^{-}(f)) = (\frac{\nu_{i_1}}{N_{i_1}} , \frac{\nu_{i_h}}{N_{i_h}}).$$
    Assume that $m$ is odd. Then by Theorem \ref{thm 7} one has $l_{1,t_1} = \gcd(m,n_1)$, which is odd since $m$ is odd. Hence $h = 1$ by Proposition \ref{prop 5} and it follows that
    $$\text{rlct}_0^{-}(f)=\text{rlct}_0^{+}(f)= \frac{\nu_{i_1}}{N_{i_1}}=\frac{1}{m} + \frac{1}{n_1}.$$
    \newline
Now assume that $m$ is even and fix $j \in \llbracket 1, g \rrbracket$. Then Theorem \ref{thm 7} gives
$$\gcd(s_j,l_{j,1}) = \gcd(l_{j,1},l_{j,2})=\dots=\gcd(l_{j,t_j-1},l_{j,t_j})=l_{j,t_j},$$
from which we deduce that there is an equivalence between the following assertions:

$$ \forall j \in \llbracket 1, h-1 \rrbracket , n_j \text{ is even  and } n_h \text{ is odd}$$

and

$$\forall j\in \llbracket 1, h-1 \rrbracket, \forall i\in \llbracket 1, t_j \rrbracket, l_{j,i} \text{ is even, and } \exists i \in \llbracket 1, t_h \rrbracket \text{ such that } l_{h,i} \text{ is odd }.$$

This, together with Proposition \ref{prop 5}, yields
$$h = \text{ min } \{ 1 \leq j \leq g \mid n_j \text{ is odd} \}.$$

\end{proof}

\begin{ex}
   We consider $f \in \mathbb{R}[x,y]$ whose Puiseux parametrization is given by
$$t \to (t^{100},t^{250}+t^{300}+t^{375} +t^{410}+t^{417})$$
and whose characteristic exponents are $(m;n_1,n_2,n_3,n_4) = (100;250,375,410,417)$. By performing the four Euclidean algorithms, we obtain that the multiplicity sequence is
$$(100,100,50,50,50,50,25,25,25,10,10,5,5,5,2,2,1,1)$$
and that the dual resolution graph, in which the exceptional divisors are numbered according to their order of appearance, is as follows.

    \begin{figure}[H]
        \centering
        \includegraphics[scale=0.35]{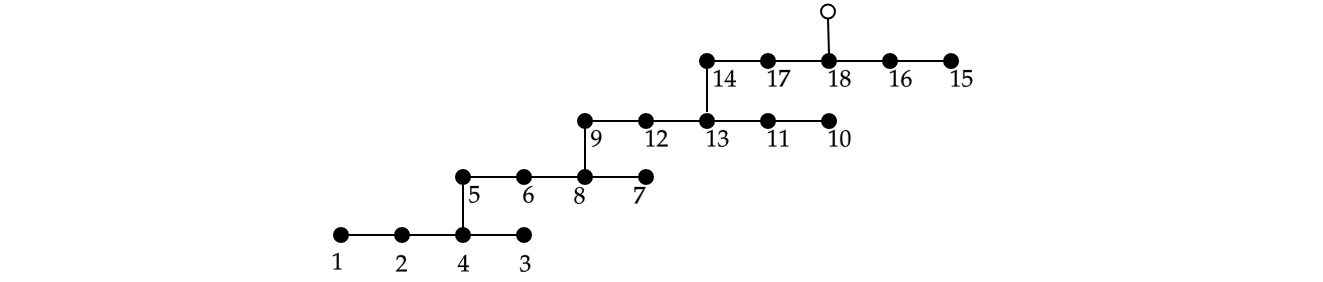}
    \end{figure}
    
   By Proposition \ref{thm 8}, the log canonical thresholds with signs are achieved at the first two rupture vertices, and we obtain
$$\frac{\nu_4}{N_4} = \frac{1}{m} +\frac{1}{n_1} = \frac{7}{500}$$
and
$$\frac{\nu_8}{N_8}=\frac{n_2+m}{e_1\overline{n}_2}= \frac{375+100}{50 \times125 + 100 \times 250}=  \frac{19}{1250}.$$

Therefore
$$\{\text{rlct}_0^{+}(f),\text{rlct}_0^{-}(f)\} = \{ \frac{7}{500}, \frac{19}{1250}\}.$$

\end{ex}

\begin{rem}
    Theorem \ref{thm 8} shows in particular that if the algebraic multiplicity of $f :(\mathbb{R}^{2},0) \to (\mathbb{R},0)$ at the origin is odd and if $f$ is analytically irreducible over $\mathbb{R}$, then one has the equality
    $$\text{rlct}_0^{+}(f)=\text{rlct}_0^{-}(f).$$
    This equality does not necessarily hold if we no longer assume that $f$ is irreducible in $\mathbb{R}\{x,y\}$. For example, consider the Newton non-degenerate and convenient germ $f(x,y) = x^{10}+x^{8}y+x^{6}y^{4}+y^{14}$ satisfying $\text{mult}_0(f) = 9$. The Newton polygon of $f$ has three compact facets $\tau_1,\tau_2,\tau_3$ and is depicted below.
    \begin{figure}[H]
        \centering
        \includegraphics[scale=0.35]{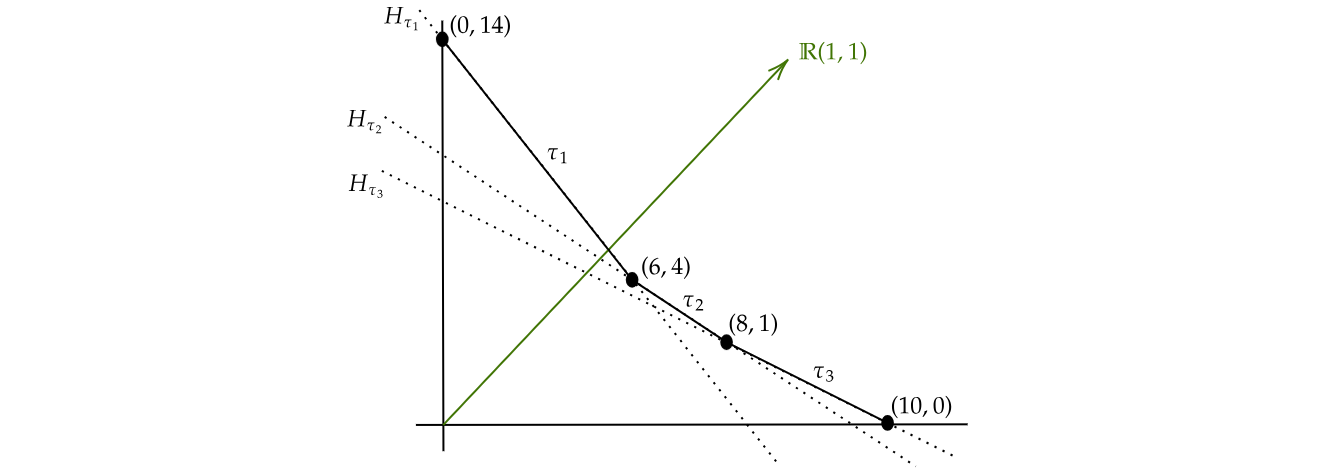}
    \end{figure}
   The supporting lines $H_{\tau_1},H_{\tau_2}$ and $H_{\tau_3}$ have equations $3y+5x=42$, $2y+3x=26$ and $2y+x=10$ respectively. One sees that the line $\mathbb{R}(1,1)$ intersects $\Gamma_c(f)$ on the facet $\tau_1$ at the point $(\frac{21}{4},\frac{21}{4})$ and that
   $$X_{\tau_1}^{+} = \{ (x,y) \in (\mathbb{R}^*)^2 \mid y^{14}+x^6y^4= 
   1 \}  \neq \emptyset$$
whereas $X_{\tau_1}^{-} = \emptyset$. Therefore $\text{rlct}_0^{+}(f) = \frac{4}{21}$. On the other hand, the line $\mathbb{R}(1,1)$ intersects $H_{\tau_2}$ at the point $(\frac{26}{5},\frac{26}{5})$ where
$$X_{\tau_2}^{-} = \{ (x,y) \in (\mathbb{R}^*)^2 \mid x^{8}y+x^6y^4= -1 \}  \neq \emptyset$$
and it follows that $\text{rlct}_0^-(f) = \frac{5}{26}$. Here, one can check that $f$ is irreducible in $\mathbb{R}[x,y]$ but that $f$ factors as a product of irreducible elements in the form $f = f_1f_2f_3$ where the $f_i$ belong to $\mathbb{R}\{x,y\}$.
 \end{rem}

\selectlanguage{english}

\bibliographystyle{amsalpha}
\bibliography{references}

\bigskip


\end{document}